\documentclass[notitlepage,reqno,12pt]{amsart}
\usepackage{latexsym,amssymb, epsfig, amsmath,amsfonts, subfigure,amsthm}

\usepackage{rotating}
\usepackage[toc,page]{appendix}
\usepackage{color}
\usepackage{multirow}
\usepackage{relsize}
\usepackage{microtype}

\usepackage[margin=1in]{geometry}
\usepackage[colorlinks, allcolors=blue]{hyperref}

\usepackage{dsfont}
\newcommand{\bneqn}{\vspace{-0.25cm}\begin{eqnarray}}
\newcommand{\eneqn}{\end{eqnarray}}

\def\rd{{\rm d}}

\numberwithin{equation}{section}
 
\newtheorem{lemma}{Lemma}[section]
\newtheorem{theorem}{Theorem}[section]

\newtheorem{prop}{Proposition}[section]
\newtheorem{remark}{Remark}[section]

\newtheorem{example}{Example}[section]

\newlength{\defbaselineskip}
\newcommand{\setlinespacing}[1]%
           {\setlength{\baselineskip}{#1 \defbaselineskip}}

\newcommand{\RR}{{\mathbb R}}
\newcommand{\ZZ}{{\mathbb Z}}

\def\E{\mathbb{E}}
\def\P{\mathbb{P}}

\newcommand{\deq}{\stackrel{\rm d}{=}}
\newcommand{\beql}[1]{\begin{equation}\label{#1}}
\newcommand{\eeq}{\end{equation}}

\newcommand{\beqal}[1]{\begin{eqnarray}\label{#1}}
\newcommand{\eeqa}{\end{eqnarray}}
\newcommand{\beq}{\begin{displaymath}}
\newcommand{\eeqno}{\end{displaymath}}
\newcommand{\bali}[1]{\begin{align}\label{#1}}
\newcommand{\eali}{\begin{align}}
\newcommand{\balino}{\begin{align*}}
\newcommand{\ealino}{\begin{align*}}

\def\tI{{\tt I}}   \def\tK{{\tt K}}  \def\tR{{\tt R}}

\newcommand{\baa}{\begin{eqnarray*}}
\newcommand{\eaa}{\end{eqnarray*}}

  \def\gam{{\gamma}}
  
\def\kap{{\kappa}} \def\lam{{\lambda}}  

 \def\vrho{\varrho}
 \def\tha{{\theta}}

\def\Llra{\Longleftrightarrow}

\def\un{\underline}
\def\une{\underline e}

\def\unw{\underline w}
\def\unx{\underline x}

 \def\unF{{\underline F}}
 \def\unH{{\underline H}}

\def\unW{{\underline W}}

\def\und0{{\underline 0}} \def\und1{{\underline 1}}

\def\cP{\mathcal P} 

\def\cR{\mathcal R} \def\cS{\mathcal S}
\def\bma{\begin{matrix}}
\def\ema{\end{matrix}}

\def\bpma{\begin{pmatrix}}
\def\epma{\end{pmatrix}}

\def\bcs{\begin{cases}} \def\ecs{\end{cases}}
\def\beac{\begin{array}{c}} \def\ena{\end{array}}
\def\beal{\begin{array}{l}} \def\beacl{\begin{array}{cl}}
\def\bealll{\begin{array}{lll}}

\def\diy{\displaystyle}

\def\beq{\begin{equation}}  
\def\eeq{\end{equation}}

\def\beal{\begin{array}{l}}
\def\ena{\end{array}}

\def\bbA{{\mathbb A}}  
 \def\bbE{{\mathbb E}}

\def\bbO{{\mathbb O}} \def\bbP{{\mathbb P}} 
\def\bbR{{\mathbb R}} \def\bbS{{\mathbb S}} 

  \def\bbZ{{\mathbb Z}}

\def\wt{\widetilde}  
\def\onwl{\operatornamewithlimits}

\newcommand{\ttl}{\Large Queueing models with random resetting}

\begin{document}

\title[]{\ttl}

\author[Dongzhou \ Huang]{Dongzhou Huang$^*$}
\address{$^*$ Department of Statistics, Colorado State University, Fort Collins, CO~~80524} \email{dongzhou.huang@colostate.edu}

					\author[Guodong \ Pang]{Guodong Pang$^{**}$}
					\address{$^{**}$Department of Computational Applied Mathematics and Operations Research,
						George R. Brown School of Engineering and Computing, 
						Rice University,
						Houston, TX 77005}
					\email{gdpang@rice.edu}

				\author[Izabella \ Stuhl]{Izabella Stuhl$^{\dag}$}
					\address{$^{\dag}$Department of Mathematics,
						Penn State University,
						University Park, PA 16802}
					\email{ius68@psu.edu}

							\author[Yuri \ Suhov]{Yuri Suhov$^{\$}$}
					\address{$^{\$}$Department of Mathematics,
						Penn State University,
						University Park, PA 16802; DPMMS, University of Cambridge, Cambridge CB3 0WB, UK; St John's College, Cambridge, Cambridge CB2 1TP, UK}
					\email{ims14@psu.edu, yms@statslab.cam.ac.uk}

\begin{abstract} 
We introduce and study some queueing models with random resetting, including Markovian and non--Markovian 
models under the first-come first-served (FCFS) discipline. The Markovian models include M/M/$r$ and M/M/1+M queues with random resetting,  
in which a continuous-time Markov chain is formulated, with transitions including a resetting to state zero in addition to arrivals and services. 
We explicitly characterize the stationary distributions of the queueing processes in these models by using parting balance equations. 
We derive expressions for standard performance measures such as the delay probability, expected queue length and waiting time, as well as probability of a customer completing service before resetting in the M/M/$r$ model, and  probability of abandonment before service or resetting in the M/M/1+M  model.

The non--Markovian models include GI/GI/1, GI/GI/$r$ and GI/GI/$\infty$ queues with random resetting to state zero at arrival times. 
For  GI/GI/1 and GI/GI/$r$ queues under the FCFS discipline, we  introduce modified 
Lindley and Kiefer--Wolfowitz recursions, respectively. Using an operator representation for these recursions, we characterize the stationary distributions via convergent series, as solutions to 
 the modified Wiener--Hopf equations. 
For  GI/GI/$\infty$ queues with resettings, we utilize a version of the
Kiefer--Wolfowitz recursion, and also characterize the corresponding stationary distribution. 
\end{abstract}


\date{\today} 

\maketitle



\allowdisplaybreaks

\section{Introduction}

In this paper we study some queueing models with random resetting, in which a queue 
clears when resettings occurs. 
Such models have many applications in service systems where machines or servers are subject to maintenance 
after some random time or periodically, or where disruptions occur due to power loss or breaking down. 
There exists a substantial literature  on how to model such phenomena, including queues in random environments, queues with disasters or catastrophes, and so on.  In this work we aim to develop a
number of queueing models with random resetting in a unified manner,  treating both Markovian and non--Markovian models.

We start with classical Markovian models with random resetting, arising from the standard 
M/M/$r$ and M/M/1+M systems.  
The queueing process in these  original models is a birth--death process with jumps $\pm 1$, whose stationary distribution is known explicitly. 
In a system with resetting, the queueing process is a continuous-time Markov chain (CTMC) with an added hop to state $0$.
Consequently, the chain with resetting is no longer a birth--death process: here the transitions 
are characterized by arrival, departure and resetting rates. 
Nevertheless, we are able to derive an explicit form for the stationary distribution.

We begin by considering the M/M/$r$ queue. For a classical M/M/$r$ queue under the first-come first-served (FCFS) discipline, the existence of a stationary distribution requires that the traffic intensity is less than one. In contrast, for the M/M/$r$ queues with random resetting under FCFS discipline, it turns out that the CTMC for the queueing process is always positive recurrent, even when the service rate equals $0$. Consequently, a stationary distribution always exists. Furthermore, it can be expressed explicitly; see Theorem~\ref{thm-MMr}. To derive the stationary distribution, we exploit the partial balance equations 
and  verify that the proposed sequence satisfies these equations.
 The resulting expression for the stationary distribution consists of two key components: (i) an analytically convenient integral representation, and (ii) a combinatorial sum, which may also be viewed as a terminating generalized hypergeometric series. Furthermore, these two components satisfy the same recursive relation but are subject to different boundary conditions; see Lemmas~\ref{lm:recursionforAn}, \ref{lm:recursionforCk}, and \ref{lm:interationAC}. Moreover, the stationary distribution exhibits a structure similar to that of the classical M/M/$r$ queues:  the distribution exhibits a Poisson-like behavior for states below $r$, while for states above $r$, it displays a geometric structure.
The M/M/1 and M/M/2 queues with resetting were studied in \cite{kumar2000transient} and \cite{kumar2002transient}, respectively, where the transient distributions were first derived and the stationary distributions were subsequently obtained as the limits of the transient distributions as time tends to infinity.
The result for the M/M/$r$ queue with a general  number of servers $r$ is novel, to the best of our knowledge.

We next study the M/M/1+M queue with random resetting under FCFS discipline, which has not been previously investigated. Following a similar approach as in the previous model, we derive an explicit expression for the stationary distribution of the queueing process; see Theorem~\ref{thm-MM1M}. In contrast to the M/M/$r$ queue with resetting, the stationary distribution in this case consists only of an analytically convenient integral representation, which takes a slight different from from that  
  appearing in the corresponding expression for the M/M/$r$ model with resetting.

 Furthermore, we observe that, when the service rate equals the abandonment rate in the M/M/1+M queue, the resulting model has the same dynamics as the M/M/$\infty$ queue with random resetting.
 Consequently, our expression for the stationary distribution reduces to the corresponding stationary distribution of the M/M/$\infty$ queue with random resetting. 
  (It is of course also the special case of M/M/$r$ queue with resetting as $r\to\infty$.) 
  Although the stationary distribution of the M/M/$\infty$ queue with random resetting was derived in \cite{bura2019transient} in a series representation using an approach similar to those in \cite{kumar2000transient} and \cite{kumar2002transient}, our result provides an alternative representation whose main component is given in terms of an analytically convenient integral. See further discussions in Remark \ref{rmk-MM1M}.

For both the M/M/$r$ and M/M/1+M queues with random resettings, we  derive the explicit expressions for the common performance measures including the delay probability, average queue length and number in the system, and expected waiting time. Moreover, for the M/M/$r$ model with resetting, we derive the 
probability that a customer eventually completes service before a resetting occurs, and for the M/M/1+M queue with resetting, we derive probability of abandonment before service or resetting.

\smallskip

Next, we analyze non--Markovian queueing models with random resetting, starting with the standard GI/GI/$1$ and
GI/GI/$r$ systems under the FCFS discipline and continuing with the GI/GI/$\infty$ system. 
For a standard GI/GI/$1$ queue, the Lindley recursion  is fundamental in studying properties of the delay/waiting-time
distribution (see Section \ref{sec-GG1-standard} for a brief review as well as  \cite[Chapter X.1]{asmussen2003}).
Similarly, for a standard GI/GI/$r$ queue, 
it is the Kiefer--Wolfowitz recursion that determines the delay/waiting--time distribution (see Section \ref{sec-GGr-standard}). We study the corresponding models where random resetting occurs at customers' arrival times, which leads to modified Lindley and Kiefer--Wolfowitz recursions (see equations \eqref{eqn-GGI-RS-W} 
and \eqref{eq:KWoR1}). 

For standard GI/GI/$1$ and GI/GI/$r$ queues under the FCFS discipline, the positive recurrence of the waiting--time process requires that the traffic intensity is less than one, and then the stationary distribution is characterized via the Wiener--Hopf equation.  
For these models 
with random resetting, we show that the corresponding modified Lindley and the Kiefer--Wolfowitz 
recursions generate positive recurrent waiting--time processes, regardless of whether their standard counterparts are positive recurrent or not. 
More importantly, the modified Lindley and the Kiefer--Wolfowitz recursions can be conveniently represented in an operator form (see equations \eqref{eqGG1:mainequation} and \eqref{OpKfr}), from which we are able to express the stationary distribution as a convergent series (see equations \eqref{eqn-FWR1} and \eqref{eq:KWoRS}-(iv)).  
As a byproduct, we give an example of a GI/GI/$1$ queue with random resetting where the interarrival times are dominated by the service times, but the waiting-time stationary distribution can be explicitly expressed, see equation \eqref{seriespsiX+}.   In this case, we are also able to derive the explicit formulas for the delay probability, and the mean and variance of the stationary waiting time (see Remark \ref{rem:pm-gg1}).

Finally, we consider GI/GI/$\infty$ models with random resetting at 
 arrival times. Here we construct a recursion for the elapsed service times for the jobs in service, by adapting 
 the Kiefer--Wolfowitz  recursion; we then use it to formulate the recursion for the  GI/GI/$\infty$ queues with resetting. 
 As a result, we derive an explicit expression for the stationary distribution (see Section \ref{sec-ISQ}).

\medskip 
{\it A review of the literature.}  The models discussed in this paper are related to several streams of the existing 
literature. First, these models are related to the stochastic clearing models studied in \cite{stidham1974stochastic,serfozo1978semi,whitt1981stationary}, where an input process (such 
as the arrival process) is intermittently and instantaneously cleared. Various clearing policies have been 
studied, e.g., clearance when the input reaches a threshold, or at i.i.d. random times independent of the 
input process.  Some of our Markovian models (the ones without service) can be regarded as stochastic clearing models 
of Poisson arrivals.
 Similarly, our non--Markovian models without 
services can be regarded as models of clearing at renewal arrival times.
However, stochastic clearing models do not generate output dynamics like our models in presence of service. 

Models considered in the present paper are also related to the queueing models with disasters or catastrophes, see, e.g, \cite{kumar2000transient,boxma2001clearing,kumar2002transient,di2003m,krishna2007transient,Yechiali2007,bohm2008note,dimou2013single,mytalas2015mx,YT17,yajima2019central,bura2019transient,sudhesh2022analysis}.  
We have already mentioned results on the M/M/1, M/M/2  and M/M/$\infty$ queues in \cite{kumar2000transient,kumar2002transient,bura2019transient}, respectively. 
In \cite{boxma2001clearing}, an M/G/$1$ queue with ``disasters" has been considered, where disasters 
occur at certain random times including (a) deterministic equidistant times, (b) random times independent of the 
queueing process, and (c) at crossings of some pre-specified level. In these works, stationary distributions of the 
workload processes have been characterized via their Laplace transforms using certain modifications of the Lindley 
recursion. 
The paper \cite{krishna2007transient} considers an M/M/$1$ queue with catastrophes, where the server 
breaks downs at i.i.d. random times, independent of the service process. At the breakdown times, all jobs are lost, and it takes
an exponential random time to repair the system. See also similar formulations of ``catastrophes" or ``clearing" in \cite{di2003m,bohm2008note,dimou2013single,mytalas2015mx,sudhesh2022analysis}. 
An associated (jump) diffusion 
approximation has been considered in \cite{di2003m,di2012double,dharmaraja2015continuous}. In 
\cite{giorno2014some}, a computational approach has been developed, for non-homogeneous Markovian single-server and 
infinite-server queueing models, whose formulation is like our Markovian models but with nonstationary transition rates. 
There are related works on M/M/1 queues with certain types of abandonment at catastrophes; see, for example, \cite{Yechiali2007, dimou2013single, YT17}. However, these models differ from the standard M/M/1+M queue with random resetting considered here. 
More general birth-death processes with catastrophes have been studied in \cite{chao2003transient,di2012double}. See also \cite{economou2003continuous} for a CTMC regulated by a point process, and \cite{stirzaker2007processes} for certain stochastic processes with random regulations.

Our study complements and provide alternative approaches to  the existing literature of Markovian queues with disasters, while   the non--Markovian models with random resetting at arrival  times are completely novel. 
We also refer the readers to some recent studies of random walks, Brownian motions and diffusions with random resetting in \cite{giorno2012reflected,stojkoski2021geometric,vinod2022time,abundo2023first,michelitsch2024random}; in some methodological aspects, they are close to our current work.

 In addition to clearing models described above, there are some recent studies of queues with random resetting in \cite{BPR22,RBP24}, 
  focused on an M/G/1 queue with service times being reset at random times whenever the service time is longer than a threshold. This concept of stochastic resetting is also exploited in random search problems, Cf. \cite{Bressloff2020,PSS2024}.

\medskip
{\it Organization of the paper.}
The paper is organized as follows. 
The Markovian models with resetting are studied first, 
with  M/M/$r$ queues in Section \ref{sec-MMr}, and M/M/1+M queues in Section \ref{sec-MM1M}. Some additional proofs for M/M/$r$ queues with resetting are collected in the Appendix,  Section \ref{sec-appendix}.
The non--Markovian models are studied next, with GI/GI/$1$ and GI/GI/$r$ queues with resettings at 
arrival times in Sections \ref{sec-GG1} and \ref{sec-GGr}, respectively, and with infinite-server queues with resettings at arrival times in Section \ref{sec-ISQ}.

\medskip

\section{The M/M/$r$ queue with random resetting} \label{sec-MMr}
\vskip .5cm

The standard assumption in Sections~\ref{sec-MMr} and \ref{sec-MM1M} is that the jobs arrive in a Poisson process at rate $\lambda  \geq 0$, the services 
times are i.i.d. exponential of rate $\mu \ge 0$, and all jobs in the system are cleared/reset after subsequent i.i.d. exponential random times of rate $\kappa >0$. The case of $\mu =0$ means that no jobs are served; in this case  we get  a stochastic clearing model of a Poisson process, Cf. \cite{stidham1974stochastic,whitt1981stationary}, and its stationary distribution is known to be geometric with parameter $\varrho := \frac{\lambda}{\lambda +\kappa}$, that is, $\pi_i=\vrho^i(1-\vrho ),\;\;i\geq 0$. 

In this section, we study an M/M/$r$ queue under the FCFS discipline with random resetting. Our analysis focuses primarily on its stationary distribution and the associated performance measures. We first describe the dynamics of the queueing system.
Let $X(t)$ be the number of jobs in the system at time $t$.  Then $\{X(t): t\ge 0\}$ is a continuous-time Markov chain (CTMC) on 
$\mathbb{Z}_+$, with the transition rates given by: 
\beq\bma i\geq 0&\to&i+1&\hbox{rate}&\lam&\hbox{(arrival),}\\
i\geq 1&\to&i-1&\hbox{rate}&(i\wedge r)\mu&\hbox{(departure),}\\
i\geq 1&\to&0&\hbox{rate}&\kap&\hbox{(resetting).}\ema\eeq
Here $x\wedge y = \min\{x,y\}$ for $x,y\in \RR$.  
Note that the CTMC $\{X(t)\}$ with $\mu>0$ is dominated by the corresponding process with $\mu=0$, and hence, it has a unique stationary distribution.



In the rest of the section, we assume $\mu >0$, and adopt the notation:
\begin{equation*}
\theta := \lambda/\mu, \quad \gamma :=\kappa/\mu.
\end{equation*}

\subsection{Stationary distribution}

This subsection is devoted to the stationary distribution. In the special case $r=2$, the authors of \cite{kumar2002transient} derived the stationary distribution by first analyzing the transient probabilities, $\pi_i(t) = P(X(t)=i)$, and then passing to the limit as $t\to\infty$. To the best of our knowledge, the case $r>2$ remains open in the literature. The corresponding general results are summarized in the following theorem.


\begin{theorem} \label{thm-MMr} Assume that $\lambda>0$, $\mu>0$, and $\kappa >0$.  
The stationary distribution $\pi$ for the {\rm{M/M/}}$r$ queue with random resetting is given by the following: 
for $k=0,1,\dots, r-1$, 
\begin{align} \label{eqn-pi-MMr-1}
\pi_{k} = \frac{ \theta^k }{k!} \left(  A_{k} -  L_{r-1,r}  C_{k} \right), 
\end{align}
and  for $k\geq r$,
\begin{align} \label{eqn-pi-MMr-2}
\pi_{k} &= \alpha^{k-r+1} \pi_{r-1} \nonumber  \\
&= \frac{\alpha^{ k-r +1}  \theta^{r-1} }{(r-1)!} \left(  A_{r-1} -  L_{r-1,r} C_{r-1} \right).
\end{align}
Here, $\theta = \lambda/\mu$, $\gamma = \kappa/\mu$,
\begin{equation}  \label{eqn-alpha-MMr}
\alpha = \frac{ r + \theta + \gamma - \sqrt{ (r + \theta + \gamma)^2  - 4\theta r } }{2r},
\end{equation}
\begin{equation}  \label{eqn-An-MMr}
A_{n} = \int_{0}^{1} \gamma (1-s)^{\gamma -1} s^n e^{-\theta s} \, ds,  \qquad \text{for $n=0,1,\dots$} ,
\end{equation}
\begin{equation} \label{eqn-Ck-MMr}
C_{k} = \sum_{l=0}^{k} \binom{k}{l} \gamma (\gamma +1) \cdots (\gamma + l-1) \theta^{-l}, \qquad \text{for $k=0,1,\dots$},
\end{equation}
and
\begin{align}  \label{eqn-L-MMr}
L_{r-1,r} 
& :=  \frac{  \theta A_{r} - r \alpha A_{r-1} }{ \theta C_{r} - r \alpha C_{r-1} }\,.
\end{align}
\end{theorem}

In view of Theorem~\ref{thm-MMr}, a natural question is whether the sequence $\{\pi_k\}$ defined in \eqref{eqn-pi-MMr-1} and \eqref{eqn-pi-MMr-2} is nonnegative and convergent. This question is addressed in Lemma~\ref{lem-pi-positive}. The proof of Lemma~2.1 relies on several auxiliary results established in Lemmas~\ref{lm:recursionforAn}--\ref{lm:interationAC}, which characterize the sequences $\{A_n\}$ and $\{C_n\}$. For brevity, the proofs of these lemmas are deferred to the Appendix.

\begin{lemma} \label{lem-pi-positive}
The sequence $\{\pi_k\}$ defined in \eqref{eqn-pi-MMr-1} and \eqref{eqn-pi-MMr-2} is always nonnegative and convergent.
\end{lemma}



\begin{lemma} \label{lm:recursionforAn}
The sequence $A_n$ defined in \eqref{eqn-An-MMr}, satisfies
\begin{equation}
(\theta + \gamma) A_0 - \theta A_1 = \gamma \label{eqinf:A0A1}
\end{equation}
and, for $n\ge 1$, 
\begin{equation} 
\theta  A_{n+1}- (n + \theta + \gamma)   A_n +  n  A_{n-1} =0. \label{eqn-An-recursion} 
\end{equation}
\end{lemma}

\begin{lemma} \label{lm:recursionforCk}
If $\theta>0$, then the sequence $C_k$ defined in \eqref{eqn-Ck-MMr}, satisfies
\begin{equation}
C_0 =1, C_1 = 1 + \frac{\gamma}{\theta}  \label{eqinf:C0C1}
\end{equation}
and, for $k \ge 1$, 
\begin{equation} 
\theta  C_{k+1}- (k + \theta + \gamma)   C_k +  k  C_{k-1} =0. \label{eqn-Ck-recursion} 
\end{equation}
Furthermore, we have
\begin{equation}
\theta C_{r} - r\alpha C_{r-1} >0.  \label{eqr:thedenominator>0}
\end{equation}
\end{lemma}

\begin{lemma}  \label{lm:interationAC}
If $\theta >0$, then for $r =1,2 \dots$, we have
\begin{equation}
\theta^{r} (  A_{r-1} C_{r} -  C_{r-1} A_{r}) = \gamma (r-1)!.   \label{eqr:keyidentityAC}
\end{equation}
\end{lemma}

With the above preparation, we are now ready to prove Theorem~\ref{thm-MMr}.

\begin{proof}[Proof of Theorem~\ref{thm-MMr}]
First, we recall that a sequence $\{\pi_k\}$ is a stationary distribution if and only if it is nonnegative and satisfies the partial balance equations (PBEs) together with the normalization condition $\sum_j \pi_j =1$, namely,
To characterize the stationary distribution, we again use the PBEs:
\begin{align}
\lambda \pi_0 &= \mu \pi_1 + \kappa \sum_{j=1}^{\infty} \pi_j, \label{eqn-MMr-BE1} \\
(\lambda + i \mu + \kappa) \pi_i &= \lambda \pi_{i-1} + (i+1) \mu \pi_{i+1}, \quad 1 \leq i < r, \label{eqn-MMr-BE2}\\
(\lambda + r \mu + \kappa) \pi_i &= \lambda \pi_{i-1} + r \mu \pi_{i+1}, \quad i \geq  r,  \label{eqn-MMr-BE3}\\
\sum_{j=0}^{\infty} \pi_j &=1, \label{eqn-MMr-BE4}\\
\pi_j &\geq 0, \quad j=0,1,2\dots. \label{eqn-MMr-BE5}
\end{align}
Therefore, it suffices to verify that the sequence $\{\pi_k\}$ given in \eqref{eqn-pi-MMr-1} and \eqref{eqn-pi-MMr-2} satisfies the above conditions. We divide the proof into two steps.
\begin{enumerate}
\item[Step 1:] We prove that the sequence $\{\pi_k\}$ satisfies the following conditions:
\begin{align}
( \theta + \gamma ) \pi_0 - \pi_1 &= \gamma,  \label{eqr:PBE1} \\
(n+1) \pi_{n+1} - (n + \theta + \gamma) \pi_n + \theta \pi_{n-1} &=0,  \qquad 1 \leq n <r,  \label{eqr:PBE2}\\
r \pi_{n+1} - (r + \theta + \gamma) \pi_n + \theta \pi_{n-1} &=0, \qquad n \geq r, \label{eqr:PBE3} \\
\sum_{n=0}^{\infty} |\pi_n| & < \infty, \label{eqr:PBE4}  \\
\pi_n &\geq 0, \quad n=0,1,2\dots.  \label{eqr:PBE5} 
\end{align}
\item[Step 2:] We establish that \eqref{eqr:PBE1}--\eqref{eqr:PBE5} imply \eqref{eqn-MMr-BE1}--\eqref{eqn-MMr-BE5}.
\end{enumerate}
Combining these two steps completes the proof.

We begin by verifying the statement in Step 2. Recalling the definitions of $\theta$ and $\gamma$, \eqref{eqn-MMr-BE2} and \eqref{eqn-MMr-BE3} are obtained by multiplying both sides of \eqref{eqr:PBE2} and \eqref{eqr:PBE3} by $\mu$, respectively. Additionally, \eqref{eqn-MMr-BE5} is identical to \eqref{eqr:PBE5}. It remains to verify \eqref{eqn-MMr-BE1} and \eqref{eqn-MMr-BE4}.

Summing \eqref{eqr:PBE2} over $n$ from $1$ to $r-1$ gives
\begin{equation*}
\sum_{n=1}^{r-1} (n+1) \pi_{n+1} - \sum_{n=1}^{r-1} (n + \theta + \gamma) \pi_n + \sum_{n=1}^{r-1} \theta \pi_{n-1} = 0,
\end{equation*}
which implies that
\begin{equation*}
r \pi_r - \theta \pi_{r-1} + \theta \pi_0 - \pi_1 - \gamma  \sum_{n=1}^{r-1} \pi_i =0.
\end{equation*}
Similarly, Summing \eqref{eqr:PBE3} over $n$ from $r$ to $\infty$, it follows after rearranging terms that
\begin{equation*}
- r \pi_r + \theta \pi_{r-1} - \gamma \sum_{n=r}^{\infty} \pi_n =0
\end{equation*}
The convergence of the series in the last display is guaranteed by \eqref{eqr:PBE4}. Combining the last two displays, we have
\begin{equation*}
(\theta + \gamma) \pi_0 - \theta \pi_1 - \gamma \sum_{n=0}^{\infty} \pi_n =0,
\end{equation*}
Together with \eqref{eqr:PBE1}, \eqref{eqn-MMr-BE4} follows. Finally, by substituting $\gamma \sum_{j=0}^{\infty}$ for $\gamma$ on the right-hand side of \eqref{eqr:PBE1} and recalling the definitions of $\theta$ and $\gamma$, we directly obtain \eqref{eqn-MMr-BE1}.

It remains to verify the statement in Step 1. By Lemma~\ref{lem-pi-positive}, \eqref{eqr:PBE4} and \eqref{eqr:PBE5} follow immediately. Moreover, substituting the expressions for $\pi_0$ and $\pi_1$ and applying Lemmas~\ref{lm:recursionforAn} and \ref{lm:recursionforCk} yields \eqref{eqr:PBE1}. It therefore remains only to prove \eqref{eqr:PBE2} and \eqref{eqr:PBE3}.

We first verify \eqref{eqr:PBE2}. The proof is divided into two cases: $1 \leq n \leq r-2$ and $n =r-1$. When $1 \leq n \leq r-2$, plugging in the expressions of $\pi_n$ in \eqref{eqn-pi-MMr-1}, it follows that
\begin{eqnarray*}
&& (n+1) \pi_{n+1} - (n + \theta + \gamma) \pi_n + \theta \pi_{n-1} \\
&=& (n+1) \frac{\theta^{n+1}}{(n+1)!} \left(  A_{n+1} - L_{r-1, r} C_{n+1} \right) - (n + \theta + \gamma) \frac{\theta^n}{n!} \left(  A_{n} - L_{r-1, r} C_{n} \right)   \\
&& + \theta  \frac{\theta^{n-1}}{(n-1)!} \left(  A_{n-1} - L_{r-1, r} C_{n-1} \right) \\
&=& \frac{\theta^n}{n!}\left( \theta A_{n+1} - (n+\theta +\gamma) A_n + nA_{n-1} \right)
- L_{r-1,r} \frac{\theta^n}{n!}\left( \theta C_{n+1} - (n+\theta +\gamma) C_n + nC_{n-1} \right) \\
&=&0,
\end{eqnarray*}
where the last equality follows from \eqref{eqn-An-recursion} and \eqref{eqn-Ck-recursion}. Similarly, when $n = r-1$,
\begin{eqnarray*}
&& r \pi_{r} - (r-1 + \theta + \gamma) \pi_{r-1} + \theta \pi_{r-2}  \\
&=& r \frac{ \alpha  \theta^{r-1} }{(r-1)!} \left(  A_{r-1} - L_{r-1, r} C_{r-1} \right) - (r-1 + \theta + \gamma) \frac{ \theta^{r-1} }{ (r-1)! } \left(  A_{r-1} - L_{r-1, r} C_{r-1} \right)  \\
&& + \theta \frac{ \theta^{r-2} }{ (r-2)! } \left(  A_{r-2} - L_{r-1, r} C_{r-2} \right) \\
&=& \frac{  \theta^{r-1} }{(r-1)!} \left( r \alpha A_{r-1} - ( r-1 +\theta +\gamma ) A_{r-1} + (r-1) A_{r-2} \right)  \\
&& - \frac{  \theta^{r-1} }{(r-1)!} L_{r-1,r} \left( r \alpha C_{r-1} - ( r-1 +\theta +\gamma ) C_{r-1} + (r-1) C_{r-2} \right).
\end{eqnarray*}
Applying \eqref{eqn-An-recursion}, it follows
\begin{eqnarray*}
 r \alpha A_{r-1} - ( r-1 +\theta +\gamma ) A_{r-1} + (r-1) A_{r-2}
 = r \alpha A_{r-1} - \theta A_{r}.
\end{eqnarray*}
Similarly,
\begin{eqnarray*}
r \alpha C_{r-1} - ( r-1 +\theta +\gamma ) C_{r-1} + (r-1) C_{r-2}
= r \alpha C_{r-1} - \theta C_{r}.
\end{eqnarray*}
Combining the last three displays, we have
\begin{eqnarray*}
&& r \pi_{r} - (r-1 + \theta + \gamma) \pi_{r-1} + \theta \pi_{r-2}  \\
&=& \frac{  \theta^{r-1} }{(r-1)!} \left(r \alpha A_{r-1} - \theta A_{r} \right)
- \frac{  \theta^{r-1} }{(r-1)!} L_{r-1,r} \left( r \alpha C_{r-1} - \theta C_{r} \right).
\end{eqnarray*}
Recalling the definition of $L_{r-1, r}$ and applying \eqref{eqr:thedenominator>0}, we obtain
\begin{equation*}
r \pi_{r} - (r-1 + \theta + \gamma) \pi_{r-1} + \theta \pi_{r-2}=0.
\end{equation*}
Combining these two cases, we establish \eqref{eqr:PBE2}.

It remains to verify \eqref{eqr:PBE3}. Plugging in the expressions of $\pi_n$ in \eqref{eqn-pi-MMr-2}, we have
\begin{eqnarray*}
&& r \pi_{n+1} - (r + \theta + \gamma) \pi_n + \theta \pi_{n-1} \\
&=& r \alpha^{n+1-r+1} \pi_{r-1} - (r + \theta + \gamma) \alpha^{n-r+1} \pi_{r-1} + \theta \alpha^{n-1-r+1} \pi_{r-1} \\
&=& \alpha^{n-r} \pi_{r-1} \left( r \alpha^2 - (r+\theta +\gamma) \alpha + \theta \right) \\
&=& 0,
\end{eqnarray*}
where the last equality is due to the definition of $\alpha$ in \eqref{eqn-alpha-MMr}.

We conclude the proof.
\end{proof}

%
%
%
%
%


\subsection{Performance measures}

In this subsection, we use the stationary distribution to investigate several performance measures. 

\subsubsection{Delay probability}
Observe that an arriving customer must wait whenever all servers are occupied. Therefore, the delay probability $P_{\mathrm{delay}}$ can be expressed as
\begin{equation*}
P_{\mathrm{delay}} = \sum_{n=r}^{\infty} \pi_n.
\end{equation*}
Since $\pi_{n}$ exhibits a geometric structure for $n \geq r$, it follows immediately that
\begin{equation*}
P_{\mathrm{delay}} = \frac{ \alpha }{1 - \alpha} \pi_{r-1}.
\end{equation*}
Furthermore, recalling the expressions for $\pi_{r-1}$ and $L_{r-1, r}$, we have
\begin{align*}
\pi_{r-1} &= \frac{ \theta^{r-1} }{(r-1)!} \left(  A_{r-1} - L_{r-1, r} C_{r-1} \right)
= \frac{ \theta^{r-1} }{(r-1)!} \left(  A_{r-1} - \frac{ \theta A_{r} - r\alpha A_{r-1} }{ \theta C_{r}  - r\alpha C_{r-1} }  C_{r-1} \right) \\
&= \frac{ \theta^{r-1} }{(r-1)!}  \frac{ \theta ( A_{r-1} C_{r} - C_{r-1} A_{r} ) }{  \theta C_{r}  - r\alpha C_{r-1} }
= \frac{ \gamma }{ \theta C_{r}  - r\alpha C_{r-1} },
\end{align*}
where the last equality follows from Lemma~\ref{lm:interationAC}. Combining the last two displays yields
\begin{equation*}
P_{\mathrm{delay}} = \frac{ \alpha }{1 - \alpha} \frac{ \gamma }{ \theta C_{r}  - r\alpha C_{r-1} }.
\end{equation*}
Although the above formula does not initially resemble the Erlang C formula for the M/M/$r$ queue, an appropriate rearrangement reveals a connection between them. Indeed, recalling the definition of $C_k$ in \eqref{eqn-Ck-MMr}, we have
\begin{eqnarray*}
&& \theta C_{r}  - r\alpha C_{r-1} \\
&=& \theta \sum_{l=0}^{r} \binom{r}{l} \gamma (\gamma +1) \cdots (\gamma + l-1) \theta^{-l}
 - r \alpha \sum_{l=0}^{r-1} \binom{r-1}{l} \gamma (\gamma +1) \cdots (\gamma + l-1) \theta^{-l} \\
&=& (\theta - r\alpha) + \theta \gamma \sum_{l=1}^{r} \binom{r}{l}  (\gamma +1) \cdots (\gamma + l-1) \theta^{-l} \\
&& - r \alpha \gamma \sum_{l=1}^{r-1} \binom{r-1}{l} (\gamma +1) \cdots (\gamma + l-1) \theta^{-l} \\
&=& (\theta - r\alpha) + (\theta - r \alpha) \gamma \sum_{l=1}^{r-1} \binom{r-1}{l} (\gamma +1) \cdots (\gamma + l-1) \theta^{-l} \\
&& + \theta \gamma \sum_{l=1}^{r} \binom{r-1}{l-1} (\gamma +1) \cdots (\gamma + l-1) \theta^{-l},
\end{eqnarray*}
where the last equality is a consequence of the identity
\begin{equation*}
\binom{r}{l} \mathds{1}_{\{ l =0,1, \dots, r \}} = \binom{r-1}{l} \mathds{1}_{\{ l =0,1, \dots, r-1 \}} + \binom{r-1}{l-1} \mathds{1}_{\{ l =1, 2, \dots, r \}}.
\end{equation*}
Therefore,
\begin{align}
P_{\mathrm{delay}} &= \bigg[ \frac{1-\alpha}{\alpha} \frac{ \theta - r \alpha }{\gamma} + \frac{1-\alpha}{\alpha} (\theta - r\alpha) \sum_{l=1}^{r-1} \binom{r-1}{l} (\gamma +1) \cdots (\gamma + l-1) \theta^{-l}  \nonumber \\
& \quad \quad \quad \quad \quad \quad \quad \quad + \frac{1-\alpha}{\alpha} \theta \sum_{l=1}^{r} \binom{r-1}{l-1} (\gamma +1) \cdots (\gamma + l-1) \theta^{-l}
\bigg]^{-1}.  \label{eqf:delayprobabilityalternative}
\end{align}
Provided that $\theta <r$, 
it can be proved that as $\gamma \rightarrow 0$ (or $\kappa \rightarrow 0$),
\begin{align*}
\frac{1-\alpha}{\alpha} \frac{ \theta - r \alpha }{\gamma} & \rightarrow 1, \\
\frac{1-\alpha}{\alpha} (\theta - r\alpha) \sum_{l=1}^{r-1} \binom{r-1}{l} (\gamma +1) \cdots (\gamma + l-1) \theta^{-l} & \rightarrow 0, \\
\frac{1-\alpha}{\alpha} \theta \sum_{l=1}^{r} \binom{r-1}{l-1} (\gamma +1) \cdots (\gamma + l-1) \theta^{-l} & \rightarrow \left( 1- \frac{\theta}{r} \right) \frac{ r! }{ \theta^r } \sum_{l=0}^{r-1} \frac{ \theta^{l} }{ l! }.
\end{align*} 
Therefore, the alternative form of $P_{\mathrm{delay}}$ in \eqref{eqf:delayprobabilityalternative} immediately shows that its limit coincides with the classical Erlang C formula for the M/M/$r$ queue.

\subsubsection{Expected queue length}
It follows from the geometric structure of $\pi_n$ for $n \geq r$ that
\begin{align*}
L_{\mathrm{queue}}:= \sum_{n=r}^{\infty} (n -r) \pi_{n} = \frac{\alpha^2}{ (1-\alpha)^2 } \pi_{r-1} = \frac{\alpha}{1-\alpha} P_{\mathrm{delay}}.
\end{align*}
This formula implies that the relationship between the expected queue length and the delay probability in the resetting case resembles that in the M/M/$r$ queue.

\subsubsection{Expected number in the system}
The expected number $L_{\mathrm{system}}$ of customers in the the M/M/$r$ queue with random resetting is defined as
\begin{equation*}
L_{\mathrm{system}} := \sum_{n=0}^{\infty} n \pi_n.
\end{equation*}
A direct calculation yields 
\begin{align*}
L_{\mathrm{system}} &= \sum_{n=0}^{r-1} n \pi_n + \sum_{n=r}^{\infty} (n-r) \pi_n + r \sum_{n=r}^{\infty} \pi_n \\
&= \sum_{n=0}^{r-1} n \pi_n + L_{\mathrm{queue}} + r P_{\mathrm{delay}}.
\end{align*}
Using a generating function approach, we can derive that
\begin{align*}
\sum_{n=0}^{r-1} n \pi_n &= \frac{\theta}{\gamma+1} - \left(  \frac{ \theta (\theta - r\alpha) }{\gamma(\gamma+1)} + \frac{ r  \theta - r(r-1)\alpha  }{\gamma+1} \right) \pi_{r-1}  \\
&= \frac{\theta}{\gamma+1} - \frac{1-\alpha}{\alpha} \left(  \frac{ \theta (\theta - r\alpha) }{\gamma(\gamma+1)} + \frac{ r  \theta - r(r-1)\alpha  }{\gamma+1} \right) P_{\mathrm{delay}}.
\end{align*}
Therefore,
\begin{equation*}
L_{\mathrm{system}} =  L_{\mathrm{queue}} + \frac{\theta}{\gamma+1} + \left[ r- \frac{1-\alpha}{\alpha} \left(  \frac{ \theta (\theta - r\alpha) }{\gamma(\gamma+1)} + \frac{ r  \theta - r(r-1)\alpha  }{\gamma+1} \right) \right]  P_{\mathrm{delay}}.
\end{equation*}
Unlike the classical M/M/$r$ queue, the formula for the expected number of customers in the system under the resetting mechanism is substantially more complex.

\subsubsection{Expected waiting time}
Let $W$ be the  waiting time of a customer entering the system in stationarity.
By Little's law, we have 
\[
\E[W]= \frac{1}{\lambda} L_{\mathrm{queue}}=  \frac{\alpha}{\lambda (1-\alpha)} P_{\mathrm{delay}}.
\]

\subsubsection{Probability that a customer eventually completes service before a resetting occurs}
Let $X$ denote the stationary number of customers in the system. We now derive the probability that a customer completes service before a resetting event occurs, conditional on the value of $X$. Denote this probability by $P_{\mathrm{CSBR}}$.

If $X = n \leq r-1$, an arriving customer is served immediately. Hence,
\begin{equation*}
\mathbb{P}\left( \text{completing service} \mid X=n \right) = \frac{ \mu }{ \mu + \kappa} = \frac{1}{1 + \gamma}.
\end{equation*}
If $X = n \geq r$, then an arriving customer must wait for the preceding $ n-r+1 $ customers to complete service. Therefore,
\begin{equation*}
\mathbb{P}\left( \text{completing service} \mid X=n \right) = \left( \frac{r \mu}{ r\mu + \kappa } \right)^{n -r +1} \frac{\mu}{\mu+\kappa} = \left( \frac{r}{r+\gamma} \right)^{n-r+1} \frac{1}{1 + \gamma}.
\end{equation*}
Combining these two cases, it follows that
\begin{align*}
P_{\mathrm{CSBR}} &= \sum_{n=0}^{r-1} \pi_{n} \frac{1}{1 + \gamma} + \sum_{n=r}^{\infty} \pi_n \left( \frac{r}{r+\gamma} \right)^{n-r+1} \frac{1}{1 + \gamma} \\
&= \frac{1}{1+\gamma} \left(1- P_{\mathrm{delay}}\right) + \frac{1}{1+\gamma} \sum_{n=r}^{\infty} \alpha^{n-r+1} \left( \frac{r}{r+\gamma} \right)^{n-r+1} \pi_{r-1} \\
&= \frac{1}{1+\gamma} \left(1- P_{\mathrm{delay}}\right) + \frac{1}{1 +\gamma} \frac{ \alpha \gamma }{ (1-\alpha)r +\gamma } \pi_{r-1} \\
&= \frac{1}{1+\gamma} \left(1- P_{\mathrm{delay}}\right) + \frac{1}{1 +\gamma} \frac{ (1-\alpha) \gamma }{ (1-\alpha)r +\gamma } P_{\mathrm{delay}} \\
&= \frac{1}{1 +\gamma} - \frac{ (1-\alpha)r + \alpha\gamma }{(1+\gamma) \left( (1-\alpha) r + \gamma \right)} P_{\mathrm{delay}},
\end{align*}
where the second equality follows by substituting the expression for $\pi_n$ from \eqref{eqn-pi-MMr-2}.


\medskip

\section{M/M/1+M queues with random resetting} \label{sec-MM1M}
\medskip

In this section, we consider the M/M/1+M queue with random resetting under the FCFS discipline, where the patience times are  i.i.d. exponentially distributed of rate $\nu$, independent from the arrival, service and and resetting processes.
Let $X(t)$ denote the number of jobs in the system at time $t \ge 0$. Then the process $\{X(t): t \ge 0\}$ is a continuous-time Markov chain taking values in $\mathbb{Z}_+ = \{0,1,2,\dots\}$ with the following transition rates:
\begin{align*}
& i \ge 0 \quad \to \quad i+1, \quad \text{rate} \quad \lambda \quad \text {(arrival)}, \\
&  i \ge 1 \quad \to \quad i-1, \quad \text{rate} \quad \mu + (i-1) \nu  \quad \text {(departure)}, \\
& i \ge 1 \quad \to \quad 0, \quad \quad\,\,\,  \text{rate} \quad \kappa \quad \text {(resetting)}. 
\end{align*}

When $\nu = 0$, the process $\{X(t): t \ge 0\}$ reduces to that in an M/M/1 queue with random resetting under the FCFS discipline. 
Owing to the dominance, we shall obtain the following property: For any $\mu\geq 0$, $\lambda \geq 0$, and $\kap >0$, 
the {\rm{CTMC}} $\{X(t): t\ge 0\}$ \ is positive recurrent and has a unique stationary distribution.  
The stationary distribution $\pi= \{\pi_i: i\in \ZZ_+\}$ has been derived in 
\cite{kumar2000transient}, by first studying the transient distribution $\pi_i(t) =P(X(t)=i)$ and then letting $t\to\infty$. 
Using the PBEs and following the same approach, the result can be readily recovered. For completeness, we state the result below.
\begin{theorem} \label{thm-MM1} 
Assume that $\nu=0$, and either $\mu > 0$ and $\kappa>0$, or $\kappa=0$ and $\mu > \lambda \ge 0$.
The  stationary distribution $\pi$ of the M/M/1 queue with resetting under the FCFS discipline is given by
\[
\pi_i = \rho^i (1-\rho), \quad i=0,1,2, \dots, 
\]
where  $\rho\in (0,1)$ is given by
\beq\label{eqn-rho}
\rho= \frac{1}{2} \Big( \lambda/\mu +\kappa/\mu +1 - \sqrt{ (\lambda/\mu +\kappa/\mu+1)^2 - 4 \lambda/\mu} \Big)\,. 
\eeq

\end{theorem}

%

We proceed to study the M/M/1+M queue with random resetting. It suffices to focus on the case where $\nu >0$. Additionally, we also assume that $\mu >0$ and $\kappa >0$.
Note that the process ${X(t)}$ with $\nu > 0$ is dominated by the process with $\nu = 0$. Based on this relation, we can conclude that the stationary distribution of ${X(t)}$ exists and is unique. Furthermore, when $\mu = \nu$, the transition rates of ${X(t)}$ coincide with those of an M/M/$\infty$ queue, implying that the dynamics of the two systems are identical. Consequently, when $\mu=\nu$, the stationary distribution of the M/M/1+M queue with random resetting shall coincide with the stationary distribution of the M/M/$\infty$ queue with resetting derived in \cite{bura2019transient}. We provide further details in the following subsection.


Before presenting our main result for the stationary distribution of the M/M/1+M queue with resetting, we pause to define necessary notation. Define
\begin{equation} \label{eqn-para-MM1M}
\theta := \frac{\lambda}{\nu}, \quad \gamma := \frac{\kappa}{\nu}, \quad \text{and} \quad \eta := \frac{\mu}{\nu}.
\end{equation}

\subsection{Stationary distribution}

\begin{theorem} \label{thm-MM1M}
Assume that $\mu, \nu, \kappa >0$. The stationary distribution $\pi$ for the M/M/1+M queue with random resetting under the FCFS discipline is given by the following: for $k =0,1,2, \dots$,
\begin{equation}
\pi_{k} = \frac{ \theta^{k} }{\prod_{i=0}^{k-1} (i + \eta) } \times  \frac{ \gamma A_{k} }{ (\gamma + \theta) A_0 - \theta A_1 },  \label{eq+M:pidef}
\end{equation}
where
\begin{equation*}
A_{k} = \int_{0}^{1} \gamma (1-t)^{\gamma -1} t^{\eta -1 +k} e^{-\theta t} dt, \qquad \text{for $k =0, 1,2, \dots$}.
\end{equation*}
\end{theorem}

\begin{proof}
To characterize the stationary distribution, we have the partial balance equations (PBEs):
\begin{align*}
\lambda\pi_0 &=\mu\pi_1+\kappa\sum\limits_{j\geq 1}\pi_j, \\
(\lambda +\mu+(k-1)\nu+\kappa )\pi_{k}&=\lambda \pi_{k-1}+(\mu + k \nu) \pi_{k+1}  , \qquad \text{for $k = 1,2, \dots$} 
\end{align*}
Using a similar argument as in the proof of Theorem \ref{thm-MMr}, we obtain that a sequence $\pi$ is the stationary distribution if it is nonnegative and convergent, and satisfies the following equations:
\begin{align}
(\theta +\gamma) \pi_{0} - \eta \pi_{1} &= \gamma,  \label{eq+M:PBE3} \\
(\theta + \eta +(k-1)+ \gamma )\pi_{k}&=\theta \pi_{k-1}+( \eta + k ) \pi_{k+1}  , \qquad \text{for $k = 1,2, \dots$}. \label{eq+M:PBE4}
\end{align}


Noting that 
\begin{align*}
A_{1} &= \int_{0}^{1} \gamma (1-t)^{\gamma -1} t^{\eta } e^{-\theta t} dt = \int_{0}^{1} \gamma (1-t)^{\gamma -1} t^{\eta -1 } \times t e^{-\theta t} dt \\
&\leq \int_{0}^{1} \gamma (1-t)^{\gamma -1} t^{\eta -1 } e^{-\theta t} dt = A_0,
\end{align*}
we have $  (\gamma + \theta) A_0 - \theta A_1 >0$. Therefore, the $\pi_k$ defined in \eqref{eq+M:pidef} are always positive. Furthermore, since $\prod_{i=0}^{k-1}  (i +\eta) \geq \eta (k-1)!$ for $k \geq 1$ and $0 \leq A_k \leq 1$ for $k \geq 1$, it is straightforward to prove that the sequence $\{\pi_{k}\}_{k=0}^{\infty}$ is convergent.
Hence, it suffices to verify that the $\pi_k$ defined in \eqref{eq+M:pidef} satisfy Equations \eqref{eq+M:PBE3} and \eqref{eq+M:PBE4}. We begin with \eqref{eq+M:PBE3}. It is straightforward that
\begin{align*}
(\theta +\gamma) \pi_{0} - \eta \pi_{1} &= (\theta +\gamma) \frac{\gamma A_0}{ (\gamma +\theta) A_0 - \theta A_1 } - \eta  \times \frac{\theta}{\eta} \frac{\gamma A_1}{ (\gamma +\theta) A_0 - \theta A_1 }  \\
&= \frac{ (\theta +\gamma) \gamma A_0 - \theta \gamma A_1 }{ (\gamma +\theta) A_0 - \theta A_1 } = \gamma.
\end{align*}
Thus, we conclude \eqref{eq+M:PBE3}.

Before the verification of \eqref{eq+M:PBE4}, we present a recursion for $A_k$, that is, for $k=1,2, \dots$,
\begin{equation*}
\theta A_{k+1} - ( \eta -1 +k +\theta + \gamma ) A_k + (\eta -1 +k) A_{k-1} =0.
\end{equation*}
In fact, using integration by parts,
\begin{align*}
A_{k} =& \int_{0}^{1} \gamma (1-t)^{\gamma -1} t^{\eta -1 +k} e^{-\theta t} dt = \int_{0}^{1} t^{\eta -1 +k} e^{-\theta t} \, d\left( - (1-t)^{\gamma} \right)  = \int_{0}^{1} (1-t)^{\gamma} \, d\left(   t^{\eta -1 +k} e^{-\theta t}\right) \\
=& (\eta -1 +k) \int_{0}^{1}  (1-t)^{\gamma} t^{\eta -1 +k -1} e^{-\theta t} dt - \theta \int_{0}^{1}  (1-t)^{\gamma} t^{\eta -1 +k } e^{-\theta t} dt \\
=&  (\eta -1 +k) \int_{0}^{1}  (1-t)^{\gamma-1} (1-t) t^{\eta -1 +k -1} e^{-\theta t} dt - \theta \int_{0}^{1}  (1-t)^{\gamma-1} (1-t) t^{\eta -1 +k } e^{-\theta t} dt \\
=& (\eta -1 +k) \int_{0}^{1}  (1-t)^{\gamma-1}  t^{\eta -1 +k -1} e^{-\theta t} dt - (\eta -1 +k) \int_{0}^{1}  (1-t)^{\gamma-1}  t^{\eta -1 +k } e^{-\theta t} dt \\
& - \theta \int_{0}^{1}  (1-t)^{\gamma-1}  t^{\eta -1 +k } e^{-\theta t} dt  + \theta \int_{0}^{1}  (1-t)^{\gamma-1}  t^{\eta -1 +k +1 } e^{-\theta t} dt \\
=& \frac{\eta -1 +k}{\gamma} A_{k-1} -  \frac{\eta -1 +k}{\gamma} A_{k} - \frac{\theta}{\gamma} A_{k} + \frac{\theta}{\gamma} A_{k+1},
\end{align*}
which is equivalent to the aforementioned recursion for $A_k$. We proceed to verify \eqref{eq+M:PBE4}. For $k=1,2,\dots$,
\begin{eqnarray*}
&& (\theta + \eta +(k-1)+ \gamma )\pi_{k}  \\
&=& (\theta + \eta +(k-1)+ \gamma ) \frac{ \theta^{k} }{\prod_{i=0}^{k-1} (i + \eta) } \frac{ \gamma A_{k} }{ (\gamma + \theta) A_0 - \theta A_1 } \\
&=& \frac{ \theta^{k} }{\prod_{i=0}^{k-1} (i + \eta) } \frac{ \gamma (\theta + \eta +(k-1)+ \gamma ) A_{k} }{ (\gamma + \theta) A_0 - \theta A_1 }  \\
&=& \frac{ \theta^{k} }{\prod_{i=0}^{k-1} (i + \eta) }  \frac{ \gamma \left(  \theta A_{k+1} + (\eta-1 +k) A_{k-1} \right) }{ (\gamma + \theta) A_0 - \theta A_1 } \\
&=& \frac{ \theta^{k+1} }{\prod_{i=0}^{k-1} (i + \eta) } \frac{ \gamma A_{k+1} }{ (\gamma + \theta) A_0 - \theta A_1 } + \frac{ \theta^{k} (\eta -1 +k) }{\prod_{i=0}^{k-1} (i + \eta) } \frac{ \gamma A_{k-1} }{ (\gamma + \theta) A_0 - \theta A_1 } \\
&=& (\eta +k) \frac{ \theta^{k+1} }{\prod_{i=0}^{k} (i + \eta) } \frac{ \gamma A_{k+1} }{ (\gamma + \theta) A_0 - \theta A_1 } + \theta \frac{ \theta^{k-1}  }{\prod_{i=0}^{k-2} (i + \eta) } \frac{ \gamma A_{k-1} }{ (\gamma + \theta) A_0 - \theta A_1 } \\
&=& (\eta +k) \pi_{k+1} + \theta \pi_{k-1},
\end{eqnarray*}
which gives \eqref{eq+M:PBE4}. We conclude the proof.
\end{proof}

\begin{remark} \label{rmk-MM1M}
\rm 
In this remark, we show that, when $\mu=\nu$, the stationary distribution of the M/M/1+M queue with random resetting coincides with that of the M/M/$\infty$ queue with resetting.
When $\mu = \nu$, it follows that $\eta = 1$. The quantities $A_n$ defined in Theorem~\ref{thm-MM1M} are identical to those in Theorem~\ref{thm-MMr}. Hence, by applying Lemma~\ref{lm:recursionforAn}, the results on the stationary distribution in Theorem~\ref{thm-MM1M} reduce to
\begin{equation*}
\pi_{k} = \frac{\theta^k}{k!} \int_{0}^{1} \gamma (1- t)^{\gamma-1} t^k e^{-\theta t} dt, \quad k =0,1,2, \dots.
\end{equation*}
This result is consistent with the stationary distribution obtained for the M/M/$\infty$ queue with resetting in \cite{bura2019transient}. 
Indeed, adopting our notation, the result in \cite[Theorem 5.1]{bura2019transient} becomes: for $k \ge 0$, 
\begin{equation*}
\pi_{k} = \sum_{j=0}^{\infty} \theta^{k+j} (-1)^{j} \frac{ (k+j)! }{j! k!} \sum_{l=0}^{k+j} \frac{(-1)^{l}}{l!  (k+j-l)!} \frac{\gamma}{\gamma +l}.
\end{equation*}
It suffices to show that the above series expressions coincide with our integral expressions. Note that, for $\gamma >0$
\begin{eqnarray*}
&& \sum_{l=0}^{k+j} \frac{(-1)^{l}(k+j)! }{l!  (k+j-l)!} \frac{\gamma}{\gamma +l}
= \sum_{l=0}^{k+j} \binom{k+j}{l} (-1)^{l} \gamma \int_{0}^{1} (1-t)^{\gamma-1+l} \, dt  \\
&=& \int_{0}^{1} \gamma (1-t)^{\gamma -1} \sum_{l=0}^{k+j} \binom{k+j}{l} (t-1)^{l} \,dt
= \int_{0}^{1} \gamma (1-t)^{\gamma -1} t^{k+j} \, dt.
\end{eqnarray*}
Combining the last two displays, it follows by Fubini's theorem that
\begin{align*}
\pi_{k} &= \sum_{j=0}^{\infty} \theta^{k+j} (-1)^{j} \frac{ 1 }{j! k!} \int_{0}^{1} \gamma (1-t)^{\gamma -1} t^{k+j} \, dt \\
&=\int_{0}^{1} \gamma (1-t)^{\gamma -1} \left( \sum_{j=0}^{\infty} \frac{ (-1)^{j} (\theta t)^{j} (\theta t)^{k} }{j!  k!} \right) \,dt  \\
&= \int_{0}^{1} \gamma (1-t)^{\gamma -1} \frac{ (\theta t)^{k} e^{-\theta t} }{k!} \, dt = \frac{\theta^k}{k!} \int_{0}^{1} \gamma (1- t)^{\gamma-1} t^k e^{-\theta t} dt. 
\end{align*}
This completes the proof of the statement.

Furthermore, while the stationary distribution of the classical M/M/$\infty$ queue follows a Poisson distribution, introducing resetting leads to a stationary distribution that remains ``Poisson-like'' in a general sense. In fact, it follows a mixed Poisson distribution. Specifically, suppose that $X|Y \sim \mathrm{Poisson}(\theta Y)$ and $Y$ has the density $\gamma (1-y)^{\gamma -1}$ for $0 < y <1$. Then, the stationary distribution $\pi$ for the M/M/$\infty$ queue with resetting satisfies $
  \pi_n = P(X =n)$.
Alternatively, the stationary distribution $\pi$ of the M/M/$\infty$ queue with resetting can be interpreted from a different perspective. Suppose that $\gamma >0$ and $\theta >0$,  define
\begin{equation*}
C = C_{\gamma, \theta} := \int_{0}^{1} (1-y)^{\gamma  -1} e^{-\theta y} \, dy,
\end{equation*}
and assume that $Y$ is a random variable with the pdf 
\begin{equation*}
  \frac{1}{C}  (1-y)^{\gamma  -1} e^{-\theta y}, \qquad 0<y<1.
\end{equation*}
Then, it is straightforward to verify
\begin{equation*}
\pi_n = C \gamma \frac{ \theta^n E\left[ Y^n \right] }{n!}.
\end{equation*}
In the special case where $\gamma=1$ (i.e., $\kappa = \mu$) and $\theta>0$, we have a more explicit formula for the stationary distribution: $\pi_n= \frac{1}{\theta} P( N >n )$, where  $N$ is a Poisson random variable with parameter $\theta$. 


\hfill $\Box$
\end{remark}

\subsection{Performance measures}
In this subsection, we investigate several performance measures for the M/M/1+M queue with random resetting, including the delay probability (namely, the probability that an arriving customer must wait for service), the expected number of customers in the system, the expected queue length, the expected actual waiting time, the expected virtual waiting time, and the probability that a customer abandons the system before receiving service or before a resetting event occurs.

\subsubsection{Delay probability}
Since there is only one server,
\begin{equation*}
P_{\mathrm{delay}} = \sum_{k=1}^{\infty} \pi_{k} = 1 - \pi_0 = 1 - \frac{ \gamma A_0 }{ (\gamma +\theta)A_0 - \theta A_1 }.
\end{equation*}

\subsubsection{Expected number in the system} \label{subsubsec:expectednumberMM1+M}
Using the explicit expression for the stationary distribution, we obtain
\begin{align*}
L_{\mathrm{system}} &:= \sum_{k=0}^{\infty} k \pi_k = \sum_{k=1}^{\infty} k \pi_k
= \sum_{k=1}^{\infty}( k -1 +\eta) \pi_k - (\eta-1)\sum_{k=1}^{\infty}  \pi_k  \\
& = \sum_{k=1}^{\infty} ( k -1 +\eta) \frac{ \theta^{k} }{\prod_{i=0}^{k-1} (i + \eta) } \frac{ \gamma }{ (\gamma +\theta)A_0 - \theta A_1 } \int_{0}^{1} \gamma (1-t)^{\gamma -1} t^{\eta -1 +k} e^{-\theta t} dt \\
& \qquad   - (\eta-1)P_{\mathrm{delay}} \\
& = \frac{ \gamma   \theta^{1-\eta} }{ (\gamma +\theta)A_0 - \theta A_1 } \int_{0}^{1} \gamma (1-t)^{\gamma -1} \sum_{k=1}^{\infty} \frac{(\theta t)^{\eta -1 +k} }{ \prod_{i=0}^{k-2} (i + \eta)  } e^{-\theta t} dt - (\eta-1)P_{\mathrm{delay}}.
\end{align*}
A routine calculation gives that
\begin{equation*}
\sum_{k=1}^{\infty} \frac{(\theta t)^{\eta -1 +k} }{ \prod_{i=0}^{k-2} (i + \eta)  }
= ( \theta t )^{\eta} + \theta^{\eta+1} t e^{\theta t} \int_{0}^{t} s^{\eta -1} e^{-\theta s} ds.
\end{equation*}
Therefore,
\begin{align*}
L_{\mathrm{system}} &= \frac{ \gamma   \theta^{1-\eta} }{ (\gamma +\theta)A_0 - \theta A_1 } \int_{0}^{1} \gamma (1-t)^{\gamma -1} ( \theta t )^{\eta} e^{-\theta t} dt  \\
& \quad \quad + \frac{ \gamma   \theta^{1-\eta} }{ (\gamma +\theta)A_0 - \theta A_1 } \int_{0}^{1} \gamma (1-t)^{\gamma -1} \theta^{\eta+1} t  \left(\int_{0}^{t} s^{\eta -1} e^{-\theta s} ds\right) dt  - (\eta-1)P_{\mathrm{delay}}  \\
&= \frac{ \gamma   \theta }{ (\gamma +\theta)A_0 - \theta A_1 } A_1 - (\eta-1) P_{\mathrm{delay}}  \\
& \quad \quad +  \frac{ \gamma   \theta^2 }{ (\gamma +\theta)A_0 - \theta A_1 }\left[  \frac{1}{\gamma(\gamma+1)} (A_0 - 2A_1 +A_2) + \frac{1}{\gamma} (A_1 -A_2) \right],
\end{align*}
where the last equality follows by applying integration by parts twice. Substituting the expression for $P_{\mathrm{delay}}$, and invoking the relation $\theta A_2 - (\eta + \theta + \gamma) A_1 + \eta A_0 =0 $, it follows after rearranging terms that
\begin{align*}
L_{\mathrm{system}} &= \frac{ \gamma (\eta -1) A_0 }{ (\gamma+1) \left[(\gamma +\theta)A_0 - \theta A_1\right]  } + \frac{\theta -\eta +1}{\gamma+1}  \\
&= \frac{\eta -1}{\gamma+1} \pi_0 + \frac{ \theta - \eta +1  }{\gamma +1}.
\end{align*}

\subsubsection{Expected queue length}
Since one customer may be in service,
\begin{align*}
L_{\mathrm{queue}} &:= \sum_{k=1}^{\infty} (k-1) \pi_k = \sum_{k=1}^{\infty} k \pi_{k} -\sum_{k=1}^{\infty} \pi_{k} = L_{\mathrm{system}} - P_{\mathrm{delay}} \\
&= \frac{ \gamma (\eta -1) A_0 }{ (\gamma+1) \left[(\gamma +\theta)A_0 - \theta A_1\right]  } + \frac{\theta -\eta +1}{\gamma+1} -1 + \frac{ \gamma A_0 }{ (\gamma +\theta)A_0 - \theta A_1 }  \\
&= \frac{ \gamma (\gamma + \eta) A_0 }{ (\gamma+1) \left[(\gamma +\theta)A_0 - \theta A_1\right]  } + \frac{\theta -\eta - \gamma}{\gamma+1}  \\
&= \frac{\gamma +\eta}{\gamma+1} \pi_0 + \frac{\theta -\eta - \gamma}{\gamma+1}.
\end{align*}

\subsubsection{Expected actual waiting time}
By Little's law, the expected actual waiting time  is given by
\[
\E[W]= \frac{1}{\lambda} L_{\mathrm{queue}} = \frac{1}{\lambda} \Big( \frac{\gamma +\eta}{\gamma+1} \pi_0 + \frac{\theta -\eta - \gamma}{\gamma+1}\Big). 
\]

\subsubsection{Expected virtual waiting time}
Let $V$ be the  virtual waiting time  given that a customer does not abandon and a resetting does not occur while waiting in queue. Let $X$ denote the stationary number of customers in the system.  
Given $X=n \geq 1$,  $V$ is distributed as  $Z_1 + \cdots + Z_n$, where $Z_1,\ldots,Z_n$ are independent exponential random variables with respective rate parameters $\mu, \mu+\nu, \dots, \mu + (n-1) \nu$.
Hence, 
\begin{align}
\E[V] &= \sum_{n=1}^\infty \pi_n \E[Z_1 + \cdots + Z_n] = \sum_{n=1}^\infty \pi_n \sum_{k=1}^{n}\frac{1}{\mu + (k-1) \nu} = \sum_{k=1}^\infty \frac{1}{\mu + (k-1) \nu} \sum_{n=k}^{\infty} \pi_n   \nonumber \\
&= \frac{1}{\nu} \sum_{k=1}^\infty \frac{1}{ \eta + (k-1) } \sum_{n=k}^{\infty} \pi_n.  \label{eqMM1+M:Vvirtualwaiting}
\end{align}
Substituting the expression for $\pi_n$ from \eqref{eq+M:pidef}, it follows that
\begin{align*}
 \sum_{n=k}^{\infty} \pi_n &= \sum_{n=k}^{\infty} \frac{ \theta^{1-\eta} \gamma }{  (\gamma +\theta) A_0 - \theta A_1 }  \frac{ 1 }{ \prod_{i=0}^{k-1} (i +\eta) } \int_{0}^{1} \gamma (1-t)^{\gamma -1} \frac{ (\theta t)^{\eta -1 +n} }{ \prod_{i=k}^{n-1} (i +\eta) } e^{-\theta t} dt \\
 &= \frac{ \theta^{1-\eta} \gamma }{  (\gamma +\theta) A_0 - \theta A_1 }  \frac{ 1 }{ \prod_{i=0}^{k-1} (i +\eta) } \int_{0}^{1} \gamma (1-t)^{\gamma -1}\left( \sum_{n=k}^{\infty} \frac{ (\theta t)^{\eta -1 +n} }{ \prod_{i=k}^{n-1} (i +\eta) }\right) e^{-\theta t} dt.
\end{align*}
Applying the identity
\begin{equation*}
\sum_{n=k}^{\infty} \frac{ (\theta t)^{\eta -1 +n} }{ \prod_{i=k}^{n-1} (i +\eta) } = (\eta + k -1) \theta^{\eta + k -1} e^{\theta t} \int_{0}^{t}  s^{\eta + k -2} e^{-\theta s }   ds,
\end{equation*}
we obtain
\begin{align*}
 \sum_{n=k}^{\infty} \pi_n = \frac{ \theta^k \gamma }{ (\gamma +\theta) A_0 - \theta A_1 } \frac{ \eta + k -1 }{ \prod_{i=0}^{k-1} (i +\eta)  } \int_{0}^{1} \gamma (1-t)^{\gamma -1} \int_{0}^{t} s^{ \eta + k -2 } e^{-\theta s} ds dt.
\end{align*}
Combining the last display with \eqref{eqMM1+M:Vvirtualwaiting}, it follows
\begin{align*}
\mathbb{E}\left[V\right] &= \frac{1}{\nu} \sum_{k=1}^{\infty}  \frac{ \theta^k \gamma }{ (\gamma +\theta) A_0 - \theta A_1 } \frac{ 1 }{ \prod_{i=0}^{k-1} (i +\eta)  } \int_{0}^{1} \gamma (1-t)^{\gamma -1} \int_{0}^{t} s^{ \eta + k -2 } e^{-\theta s} ds dt  \\
&= \frac{1}{\nu} \frac{ \theta^{1 -\eta} \gamma }{ (\gamma +\theta) A_0 - \theta A_1 } \int_{0}^{1} \gamma (1-t)^{\gamma -1} \int_{0}^{t} s^{ -1 } e^{-\theta s} \left( \sum_{k=1}^{\infty} \frac{ (\theta s)^{\eta + k -1 } }{\prod_{i=0}^{k-1} (i +\eta)} \right) ds dt.
\end{align*}
As before, by expressing the infinite series in integral form, it follows after rearranging terms that
\begin{equation*}
\mathbb{E}\left[V\right] = \frac{1}{\nu} \frac{ \theta \gamma }{ (\gamma +\theta) A_0 - \theta A_1 } \int_{0}^{1} \gamma (1-t)^{\gamma -1} \int_{0}^{t} s^{ -1 } \int_{0}^{s} u^{\eta -1} e^{-\theta u} du ds dt.
\end{equation*}

\subsubsection{Probability of abandonment before service or resetting}
Let $A$ and $R$ denote the patience time of an arriving customer and the time until the occurrence of a random resetting event, respectively. It is obvious that $A \sim \mathrm{Exp}(\nu)$ and $R \sim \mathrm{Exp}(\kappa)$. 

We aim to evaluate 
\begin{equation*}
P_{\mathrm{abandon}} := \mathbb{P} \left( A < V \wedge R \right),
\end{equation*}
where $V \wedge R=\min\{V,R\}$.

If $X=0$, it is immediate that
\begin{equation*}
\mathbb{P}\left( A < V \wedge R \mid X=0 \right) =0.
\end{equation*}
If $X= n \geq 1$, then
\begin{align*}
\mathbb{P}\left( A < V \wedge R \mid X=n \right) &= P\left( A <  R \mid X=n \right) - P\left( V < A < R \mid X=n \right)  \\
&= \frac{\nu}{ \kappa + \nu } - P\left( V < A < R \mid X=n \right).
\end{align*}
By exploiting the distributions of $A$, $R$, and $S$, it follows that 
\begin{eqnarray*}
&& \mathbb{P}\left( V < A < R \mid X=n \right) \\
&=& \int_{0}^{\infty} \cdots \int_{0}^{\infty} \int_{z_1 + \cdots +z_n}^{\infty} \left( e^{-\nu(z_1 + \cdots +z_n)} - e^{-\nu r}\right) \kappa e^{-\kappa r} dr  \\
&& \qquad \qquad \times\mu e^{-\mu z_1} (\mu + \nu)e^{-(\mu+\nu)z_2} \cdots (\mu + (n-1)\nu)e^{-(\mu+(n-1) \nu)z_n} dz_1 \cdots dz_n  \\
&=& \int_{0}^{\infty} \cdots \int_{0}^{\infty} \frac{\nu}{\kappa+\nu} e^{-(\kappa+\nu) ( z_1 + \cdots +z_n )}  \\
&& \qquad \qquad \times \mu e^{-\mu z_1} (\mu + \nu)e^{-(\mu+\nu)z_2} \cdots (\mu + (n-1)\nu)e^{-(\mu+(n-1) \nu)z_n} dz_1 \cdots dz_n  \\
&=& \frac{\nu}{\kappa+\nu} \frac{\mu}{\mu + \kappa+\nu} \frac{\mu + \nu}{\mu + \kappa+ 2\nu} \cdots \frac{\mu + (n-1)\nu}{\mu + \kappa+ n\nu}.
\end{eqnarray*}
Combining the last two displays, we have
\begin{align*}
\mathbb{P}\left( A < V \wedge R \mid X=n \right) &= \frac{\nu}{\kappa+\nu} - \frac{\nu}{\kappa+\nu} \frac{\mu}{\mu + \kappa+\nu} \frac{\mu + \nu}{\mu + \kappa+ 2\nu} \cdots \frac{\mu + (n-1)\nu}{\mu + \kappa+ n\nu} \\
&= \frac{1}{\gamma+1} - \frac{1}{\gamma+1} \frac{\eta}{\eta +\gamma +1} \frac{\eta+1}{\eta +\gamma +2} \cdots \frac{\eta+ n-1}{\eta +\gamma +n} \\
&= \frac{1}{\gamma+1} - \frac{1}{\gamma+1} \prod_{i=0}^{n-1} \frac{ (\eta + i) }{(\eta +\gamma +i+1)}.
\end{align*}
Therefore,
\begin{align*}
P_{\mathrm{abandon}} &= \sum_{n=1}^{\infty} \left(  \frac{1}{\gamma+1} - \frac{1}{\gamma+1} \prod_{i=0}^{n-1} \frac{ (\eta + i) }{(\eta +\gamma +i+1)}  \right) \pi_n  \\
&= \frac{1}{\gamma +1} P_{\mathrm{delay}} - \frac{1}{\gamma +1} \sum_{n=1}^{\infty} \frac{ \theta^n }{ \prod_{i=0}^{n-1} ( \eta +\gamma +i+1 ) } \frac{\gamma A_n}{ (\gamma+\theta) A_0 - \theta A_1 }.
\end{align*}
Using a similar argument to that in Subsubsection~\ref{subsubsec:expectednumberMM1+M}, we obtain the following alternative expression for $P_{\mathrm{abandon}}$:
\begin{eqnarray*}
&& P_{\mathrm{abandon}}   \\
&=& \frac{1}{\gamma+1} P_{\mathrm{delay}} - \frac{ \gamma \theta }{ (\gamma+1) \left[ (\gamma+\theta) A_0 - \theta A_1 \right] }
\int_{0}^{1} \gamma (1-t)^{\gamma-1} t^{-\gamma -1} \left(  \int_{0}^{t} s^{\eta + \gamma} e^{-\theta s} ds \right) dt  \\
&=& \frac{1}{\gamma+1} P_{\mathrm{delay}} - \frac{ \gamma \theta }{ (\gamma+1) \left[ (\gamma+\theta) A_0 - \theta A_1 \right] } \int_{0}^{1} \int_{0}^{1} \gamma (1-t)^{\gamma-1} t^{\eta} s^{\eta +\gamma} e^{-\theta st} ds dt, 
\end{eqnarray*}
where the last equality follows by the change of variables $s\mapsto st$. 
Based on the expressions in the two displays above, it appears that the probability of abandonment before receiving service or before a resetting event occurs does not have a simple closed-form expression, in contrast to the classical M/M/1+M queue.

\vskip .8cm

\section{The GI/GI/$1$ queues with random resetting at arrival times}
 \label{sec-GG1} \vskip .5cm

In this section, we consider a GI/GI/$1$ queue under the FCFS discipline with random resetting 
at arrival times, particularly, focusing on the waiting times (delays) of jobs in the system.
We will be using \cite{asmussen2003} as a main reference book; the original works containing 
related results can be traced via comments and the bibliography in \cite{asmussen2003}.

\subsection{\it The Lindley recursion for {\rm{GI/GI/}}$1$ queue}\label{sec-GG1-standard}
The key ingredient of  the GI/GI/$1$ model is a sequence of real random variables (RVs) 
$X_n$, $n\geq 0$, where $X_n=V_n-U_n$, $U_n$ is the $n$th inter-arrival time and $V_n$
the $n$th service time. It is assumed that these RVs $X_n$ are IID, with a common  
cumulative distribution function (CDF) $F_X$. In all of Sections \ref{sec-GG1}--\ref{sec-ISQ}, we assume that
$F_X$ is a proper CDF on $\mathbb{R}$. The latter signifies that $\lim\limits_{x\to -\infty}F_X(x)=0$,
$\lim\limits_{x\to\infty}F_X(x)=1$, i.e., that the RVs $X_n$ take only finite values.  

Let $W_n$ be the $n$th waiting time. The Lindley recursive equation (originated in \cite{Lindley52}) 
states:
\beq\label{eqn-GGI-W}
W_{n+1}= (W_n + X_n)^+, \quad n \ge 0, \eeq
with some given initial RV $W_0\ge 0$, assumed to be independent of $\{X_n\}$. Here and below, 
we set $Y^+=0\vee Y$. Then 
$\{W_n,\;n\geq 0\}$
is a discrete-time Markov chain (DTMC) on $\mathbb{R}_+=[0,\infty )$.

It is known (see, e.g., \cite[Chapter X.1]{asmussen2003}) that if $\E[X]<0$, there exists a 
unique stationary distribution of DTMC $\{W_n\}$. (In fact, for $\E[X]<0$, the DTMC
$\{W_n\}$ is Harris ergodic.) The stationary distribution is characterized by a proper CDF $F_W$  on
$\mathbb{R}_+$ determined as a unique solution to the stationary Wiener--Hopf (WH) equation
\beq\label{eqn-SWH}F_W(t) = (F_X*F_W)(t) {\bf1}_{\mathbb{R}_+}(t),\quad t\in\mathbb{R}. \eeq
Here and below, $G_1*G_2$ means the convolution of CDFs $G_i$: 
$$G_1*G_2(t)=\int_\RR G_1(t-y) d G_2(y)=G_2*G_1(t),\quad t\in\RR .$$
When $\E[X]\geq 0$, \eqref{eqn-SWH} has no solution among proper CDFs (again, 
see \cite[Chapter X.1]{asmussen2003}).  

A stochastic version of equation \eqref{eqn-SWH} reads  
\beq W \deq  (W+X)^+\,.\eeq
Here $X$ and $W$ are `generic' RVs with CDFs $F_X$ and $F_W$, respectively, independent 
of each other, and $\deq$ means ``equality in distribution''.



\subsection{\it The modified Lindley recursion for a {\rm{GI/GI/}}$1$ queue with resetting.}
 \label{sec-GG1-reset} \vskip .5cm

We consider a GI/GI/$1$ model where random resettings occur independently at arrival 
times. That is, the $(n+1)$st reset waiting time $W^\tR_{n+1}$ either continues 
as in \eqref{eqn-GGI-W} with probability $q\in (0,1)$, or is set to be $0$ with probability $1-q$, independently of  
$(X_k,W^\tR_k)$ with $0\leq k\leq n$. Recursively, it can be expressed as follows:
\beq\label{eqn-GGI-RS-W}
W^\tR_{n+1}=Z_{n+1}(W^\tR_n+X_n)^+, \quad n \geq 0.
\eeq 
Here $\{Z_n: n\ge 1\}$ is a sequence of IID Bernoulli RVs with probability $\mathbb{P}
(Z_n=0)=q=1-\mathbb{P} (Z_n=1)$, independent of $\{X_n\}$. Equivalently, we can write
\beq\label{eqn-GGI-RS-W-1}
W^\tR_{n+1}=\begin{cases} 0,&\hbox{with probability $q$,}\\
(W^\tR_n+X_n)^+,&\hbox{with probability $1-q$,}\end{cases}\eeq
independently of $(X_k,W^\tR_k)$ with $0\leq k\leq n$.

Equations \eqref{eqn-GGI-RS-W} and \eqref{eqn-GGI-RS-W-1} are referred to as a modified Lindley recursion  
with resetting. 
The sequence $\{W^\tR_n\}$ forms a DTMC on $\mathbb{R}_+$.  We show that it is 
Harris ergodic. This implies that the DTMC 
$\{W^\tR_n\}$ has a unique stationary distribution, and that the corresponding CDF, denoted 
by $F_{W^\tR}$, is proper on $\mathbb{R}_+$ and satisfies $F_{W^\tR}(0)>0$. 


\begin{prop}
For any $q\in (0,1)$ and a sequence of \ {\rm{IID}} \ {\rm{RV}}s \ $\{X_n\}$,  
the {\rm{DTMC}} \ $\{W^\tR_{n}\}$ \ is Harris ergodic. 
\end{prop}

\begin{proof}
Define $T= \inf\{n \geq 0 : W_{n}^{\texttt{R}}=0 \}$, i.e., $T$ is the first hitting time of $0$. It is straightforward that for $n \geq 1$, 
\begin{align*}
\mathbb{P}( T> n) \leq \mathbb{P}\left( Z_1 \neq 0, \cdots, Z_{n} \neq 0 \right) = (1-q)^{n}.
\end{align*}
Therefore,
\begin{equation*}
\mathbb{E}\left[ T \right] = \sum_{n=0}^{\infty} \mathbb{P}( T> n) \leq \frac{1}{q} < \infty.
\end{equation*}
Hence, the time for the chain to return to state $0$ has a finite mean, and therefore the state $0$ forms a regeneration set. Then, Harris ergodicity is straightforward.

Observe that the process is regenerative (possibly after the first cycle in case the system starts
from $W^\tR_0>0$), with the cycles $(W^\tR_{1}, \dots, W^\tR_{T})$ where $W^\tR_{1}=0$ and
$W^\tR_{T}$ for $T$ being geometric of parameter $q$. That is, the state $0$ forms a regeneration set.
Hence, the time for the chain to return
to state $0$ has a finite mean. Then, Harris ergodicity is straightforward.
\end{proof}

The stationary CDF $F_{W^\tR}$ is identified as a solution to a stationary WH equation
with resetting
\beq\label{eqn-WHeqn-R}
F_{W^\tR}(t) = \Big[ q +(1-q) (F_{W^\tR}*F_X)(t) \Big] {\bf1}_{\RR_+}(t), \quad t\in \RR, 
\eeq
or to its stochastic analog
\beq W^\tR \deq Z(W^\tR+X)^+\,.\eeq
Here $X$ and $W^\tR$ are `generic' RVs with CDFs $F_X$ and $F_{W^\tR}$, respectively,
and $Z$ is a Bernoulli RV with $\mathbb{P}(Z=0)=q=1-\mathbb{P}(Z=1)$. Furthermore, the RVs $X$, 
$W^\tR$ and $Z$ are independent, and, as before,  $\deq$ means ``equality in distribution''.

\subsection{\it The operator calculus for a {\rm{GI/GI/}}$1$ queue with resetting.}
 \label{sec-GG1-opera} \vskip .5cm

To study the stationary distribution of the GI/GI/$1$ queue with resetting, it is useful to introduce the following operator concepts. We define $\mathfrak{F}$ to be the space of all proper CDFs. It is well-known that $\mathfrak{F}$ is closed under convolution and is convex with respect to addition, where addition is understood as pointwise addition of CDFs. Furthermore, we define an operator \texttt{K} from $\mathfrak{F}$ to $\mathfrak{F}$ by
\begin{equation*}
\begin{array}{rccc}
\texttt{K} : &\mathfrak{F} &\longrightarrow& \mathfrak{F} \\[6pt]
&H(x) &\mapsto& (H * F_{X})(x) \,\mathbf{1}_{\mathbb{R}_+}(x).
\end{array}
\end{equation*}
It is straightforward to show that $ (H * F_{X})(x) \,\mathbf{1}_{\mathbb{R}_+}(x)$ is a proper CDF, and therefore \texttt{K} is a valid operator on $\mathfrak{F}$. Additionally, we can prove that $\texttt{K}$ is convex-combination preserving, in the sense that
\begin{equation*}
\tK \left(\sum_{i=1}^n\alpha_iG_i\right)=\sum\limits_{i=1}^n\alpha_i\tK\left( G_i\right) ,
\end{equation*}
for any coefficients $\alpha_1,\dots ,\alpha_n\geq 0$ with $\sum_{i=1}^n\alpha_i=1$ and CDFs $G_1$, $\dots$, $G_n$. Moreover, a corresponding infinite-dimensional version also holds, namely,
\begin{equation}
\texttt{K}\left( \sum_{i=1}^{\infty}\alpha_iG_i \right) = \sum_{i=1}^{\infty} \alpha_{i} \texttt{K}\left( G_i \right),  \label{eqGG1:continuityK}
\end{equation}
for any nonnegative sequence $\{\alpha_i\}_{i=1}^{\infty}$ with $\sum_{i=1}^{\infty} \alpha_i =1 $ and any sequence $\{G_{i}\}_{i=1}^{\infty}$ of CDFs. Indeed, it follows from the definition of $\texttt{K}$ that
\begin{eqnarray*}
&&\texttt{K}\left( \sum_{i=1}^{\infty}\alpha_iG_i \right) (x) = \left[ \left(  \sum_{i=1}^{\infty}\alpha_iG_i \right) * F_{X}\right](x) \,\mathbf{1}_{\mathbb{R}_+}(x) \\
&=& \int_{-\infty}^{\infty} \left(  \sum_{i=1}^{\infty}\alpha_iG_i \right) (x -y) \, d F_{X}(y) \,\mathbf{1}_{\mathbb{R}_+}(x) 
= \int_{-\infty}^{\infty} \sum_{i=1}^{\infty}\alpha_i G_{i}(x-y)  \, d F_{X}(y) \,\mathbf{1}_{\mathbb{R}_+}(x) \\
&=&  \sum_{i=1}^{\infty} \alpha_i \int_{-\infty}^{\infty} G_{i}(x-y)  \, d F_{X}(y) \,\mathbf{1}_{\mathbb{R}_+}(x) 
= \sum_{i=1}^{\infty} \alpha_i \left(G_{i} * F_{X}\right)(x) \,\mathbf{1}_{\mathbb{R}_+}(x) \\
&=& \sum_{i=1}^{\infty} \alpha_{i} \texttt{K}\left( G_i \right)(x),
\end{eqnarray*}
where the third last equality follows by the monotone convergence theorem. Finally, we denote by $\texttt{I}$ the identity operator on $\mathfrak{F}$, that is, $\texttt{I}$ maps each proper CDF to itself. Additionally, we define $ \texttt{K}^{n} $ to be the $n$-fold composition of $\texttt{K}$; in particular,  $ \texttt{K}^{0} = \texttt{I} $.

With the above preparation, we proceed to derive the stationary distribution by solving for a proper CDF $H(x)$ satisfying
\begin{equation*}
H(t) = \Big[ q +(1-q) (H*F_X)(t) \Big] {\bf1}_{\RR_+}(t), \quad t\in \RR.  
\end{equation*}
In terms of the operator \texttt{K}, this equation is equivalent to
\begin{equation}
H = q \mathbf{1}_{\mathbb{R}_{+}} + (1-q) \texttt{K}\left( H \right). \label{eqGG1:mainequation}
\end{equation}
To solve this equation, we consider the following series 
\begin{align}
q \big(\tI -(1-q)\tK\big)^{-1}({\bf 1}_{\RR_+})
:= & q\sum_{j = 0}^{\infty}(1-q)^j\tK^j \left({\bf 1}_{\RR_+}\right) = \sum_{j = 0}^{\infty} q (1-q)^j\tK^j  \left({\bf 1}_{\RR_+}\right)  \nonumber \\
 =& q{\bf 1}_{\RR_+}+ q(1-q)F_X{\bf 1}_{\RR_+}+q(1-q)^2(F_X*(F_X{\bf 1}_{\RR_+})){\bf 1}_{\RR_+}  \nonumber \\
& \  +q(1-q)^3(F_X*((F_X*(F_X{\bf 1}_{\RR_+})){\bf 1}_{\RR_+})){\bf 1}_{\RR_+}+\cdots,  \label{eqGG1:solutionseries}
\end{align}
where in the last equality, we applied the fact that $ {\bf 1}_{\RR_+}$ is the identity element with respect to convolution. First, we note that the function defined by the above series is a proper CDF, since $\texttt{K}^{j} \left( \mathbf{1}_{ \mathbb{R}_{+} } \right) $ are proper CDFs 
for $j = 0, 1, 2,\dots$, and the corresponding coefficients are nonnegative and sum to $1$. Subsequently, it follows by \eqref{eqGG1:continuityK} that
\begin{align*}
& q \mathbf{1}_{\mathbb{R}_{+}} + (1-q) \texttt{K}\left( \sum_{j = 0}^{\infty} q (1-q)^j \tK^j  \left({\bf 1}_{\RR_+}\right)  \right) 
= q \mathbf{1}_{\mathbb{R}_{+}} + (1-q) \sum_{j = 0}^{\infty} q (1-q)^j \texttt{K} \left( \texttt{K}^{j} \left( \mathbf{1}_{ \mathbb{R}_{+} } \right) \right)  \\
=& q \mathbf{1}_{\mathbb{R}_{+}} + \sum_{j = 0}^{\infty} q (1-q)^{j+1}  \texttt{K}^{j+1} \left( \mathbf{1}_{ \mathbb{R}_{+} } \right) 
= q \mathbf{1}_{\mathbb{R}_{+}} + \sum_{j = 1}^{\infty} q (1-q)^{j}  \texttt{K}^{j} \left( \mathbf{1}_{ \mathbb{R}_{+} } \right) \\
=& \sum_{j = 0}^{\infty} q (1-q)^{j}  \texttt{K}^{j} \left( \mathbf{1}_{ \mathbb{R}_{+} } \right).
\end{align*}
In summary, the series \eqref{eqGG1:solutionseries} defines a proper CDF and provides a solution to \eqref{eqGG1:mainequation}. Furthermore, it follows by \eqref{eqn-WHeqn-R} that $F_{W^\tR}$ also satisfies \eqref{eqGG1:mainequation}. If it can be shown that \eqref{eqGG1:mainequation} admits a unique bounded solution, then we may conclude that $F_{W^\tR}$ coincides with the series \eqref{eqGG1:solutionseries}, thereby yielding an explicit expression for the CDF of the stationary distribution. The preceding deduction is summarized in the following theorem.

\begin{theorem} \label{thm-GG1} For any $q\in (0,1)$, the equation \eqref{eqGG1:mainequation} has a unique bounded solution. Since both $F_{W^{\texttt{R}}}$ and the CDF defined by the series \eqref{eqGG1:solutionseries} satisfy \eqref{eqGG1:mainequation}, they must coincide. That is, 
\begin{align} 
F_{W^{\texttt{R}}} =& \sum_{j = 0}^{\infty} q (1-q)^j\tK^j  \left({\bf 1}_{\RR_+}\right)  \nonumber  \\
 =& q{\bf 1}_{\RR_+}+ q(1-q)F_X{\bf 1}_{\RR_+}+q(1-q)^2(F_X*(F_X{\bf 1}_{\RR_+})){\bf 1}_{\RR_+}  \nonumber  \\
& \  +q(1-q)^3(F_X*((F_X*(F_X{\bf 1}_{\RR_+})){\bf 1}_{\RR_+})){\bf 1}_{\RR_+}+\cdots.\label{eqn-FWR1} 
\end{align}
\end{theorem}

\begin{proof} 
We only need to check uniqueness. Let $G_1\,:\RR\to\RR$ and  $G_2\,:\RR\to\RR$ be two bounded 
functions such that $G_i(t) = \Big[ q +(1-q) (G_i*F_X)(t) \Big] {\bf1}_{\mathbb{R}_+} (t)$, for $t\in \RR$, $i=1,2$.
Setting $G:=G_1-G_2$,  then  
$$G(t)=(1-q)(G*F_X)(t){\bf 1}_{\RR_+}(t),\quad t\in\RR.$$
To establish the uniqueness, it suffices to show that $  G \equiv 0$. To this end, define $\Gamma :=\sup_{t\in\RR} |G(t)|$. If $\Gamma >0$, then
\begin{align*}
\Gamma & = \sup_{t \in \mathbb{R}} | G(t) | = \sup_{t \in \mathbb{R}} \left| (1-q) \int_{-\infty}^{\infty} G(t -y) \, dF_{X}(y) \, \mathbf{1}_{ \mathbb{R}_{+} } (t) \right|  \\
&\leq (1-q) \sup_{t \in \mathbb{R}} \int_{-\infty}^{\infty} \left| G(t -y)  \right|\, dF_{X}(y) \, \mathbf{1}_{ \mathbb{R}_{+} } (t) \\
& \leq (1-q) \sup_{t \in \mathbb{R}} \int_{-\infty}^{\infty} \Gamma \, dF_{X}(y) \, \mathbf{1}_{ \mathbb{R}_{+} } (t) \\
& = (1-q) \Gamma.
\end{align*}
This is a contradiction, since $1-q <1$. Therefore, $\Gamma =0$, which implies $G \equiv 0$.
\end{proof}\vskip .5cm

We now turn to the case where RV $X\geq 0$. In this case, the CDF $F_{X}(x)$ is positive only if $x\geq 0$. Therefore, $ F_{X} \mathbf{1}_{   \mathbb{R}_{+} } = F_{X} $. Then, $F_{W^\tR}$ determined by the series in \eqref{eqn-FWR1} reduces to 
\beq\label{seriesRWX+}\beal F_{W^\tR}= q{\bf 1}_{\RR_+}+ q(1-q)F_X  +q(1-q)^2F_X*F_X
+q(1-q)^3F_X*F_X*F_X+\cdots.\ena\eeq

Let $\phi(\theta) = \E[e^{ {\rm i} \theta X}]$ and $\psi(\theta) = \E[e^{{\rm i} \theta W^\tR}]$ be the characteristic functions (CFs) of $X$ and $W^\tR$, respectively, for $\theta\in \RR$. Then 
\eqref{seriesRWX+} is equivalent to 
\beq\label{seriespsiX+}\beac\psi (\tha )=q\sum\limits_{k= 0}^{\infty}(1-q)^k(\phi (\tha ))^k
=\dfrac{q}{1-(1-q)\phi (\tha )}.\ena\eeq

We obtain the following assertion.    

\begin{theorem} \label{thm-GG1-Xpos}
For any $q \in (0,1)$ and {\rm{RV}} $X\ge 0$, the {\rm{RV}} $W^\tR$, corresponding to the stationary distribution of the GI/GI/$1$ queue with resetting, is a proper {\rm{RV}} with values in $[0,\infty )$. Furthermore, its {\rm{CDF}} $F_{W^\tR}$ and {\rm{CF}} $\psi$ are given by the expressions in \eqref{seriesRWX+} and \eqref{seriespsiX+}.

\end{theorem}

As an example, we next consider the case $X$ taking integer values. \vskip .5cm

\begin{example}
\rm 
 We start with the simplest scenario:  $X=1$. For simplicity, suppose that $W_0=0$. 
In the model without resetting, the waiting time $W_n=n\to\infty$ as $n\to\infty$. Cf. \eqref{eqn-GGI-W}.
In the presence of resetting, we have a DTMC  $\{W^\tR_{n}\}$ on the state space $ \ZZ_+$, with transition
probabilities $P=(P_{ij})$:  $P_{i,0} =q$ and $P_{i,i+1} =1-q$ for $i\ge 0$. 
The stationarity condition $\pi P= \pi$ gives $q \sum\limits_{i\geq 1}\pi_i = \pi_0$ and 
$(1-q) \pi_i = \pi_{i+1}$ 
for $i \ge 0$. Thus, we obtain that $\pi_i = q(1-q)^i$, $i\ge 0$, i.e., the RV $W^\tR$  is geometric.     
It is straightforward that the CDF \ $F_{W^\tR}$ and CF \ $\psi$ satisfy  \eqref{seriesRWX+} and \eqref{seriespsiX+}. \vskip .5cm
\end{example}

\begin{example} 
\rm Next, assume that $X$ takes values $k=1,2,\ldots$, with probabilities $p_k$ and
take again $W_0=0$. Then the  DTMC  $\{W^\tR_{n}\}$ has transition probabilities 
$P_{i,0}=q$ for $i\ge 0$, and $P_{i,j} = (1-q) p_{j-i}$ for $i\ge 0$ and $j\ge i+1$; the remaining 
entries equal to zero.  Let $\pi= (\pi_0,\pi_1,\dots)$ be the stationary distribution. Then we get 
$q \sum\limits_{i=0}^\infty \pi_i = \pi_0$, and $(1-q)  \sum\limits_{j=0}^{i-1} \pi_j p_{i-j} = \pi_i$ for $i \ge 1$. 
From this we obtain the equation for the characteristic functions: $\psi(\theta) = 1+ (1-q) \psi(\theta)
\phi(\theta)$, which results in \eqref{seriespsiX+}.
\end{example}

\begin{remark} \label{rem:pm-gg1}
\rm 
The series representation \eqref{eqn-FWR1} of the stationary distribution for the GI/GI/1 queue with random resetting enables us to derive an expression for the delay probability $P_{\mathrm{delay}}$ in terms of probabilities involving independent copies $X_1, X_2, \dots$ of $X$. To this end, recall that if $F_1$ and $F_2$ denote CDFs of two independent random variables $Y_1$ and $Y_2$, respectively, then $F_1 * F_2$ is the CDF of $Y_1 + Y_2$. Furthermore, $F_1 \mathbf{1}_{\mathbb{R}_+}$ is the CDF of $Y_1^+$. Combining these observations with \eqref{eqn-FWR1}, we obtain
\begin{align*}
F_{W^{\texttt{R}}} (0) &=  q +  q(1-q) \mathbb{P}(X_1 \leq 0) + q(1-q)^2 \mathbb{P} \left(  X_{1}^{+} + X_2 \leq 0 \right)  \\
& \quad + q (1-q)^3 \mathbb{P}\left(  \left(X_1^{+} + X_2\right)^{+} + X_3 \leq 0  \right) \\
& \quad + q (1-q)^4 \mathbb{P}\left(  \left(\left(X_1^{+} + X_2\right)^{+} + X_3 \right)^{+} + X_4  \leq 0  \right) + \cdots.
\end{align*}
Noting that $P_{\mathrm{delay}} = 1 - F_{W^{\texttt{R}}} (0)$, we obtain an expression for the delay probability. For the expected stationary waiting time $\mathbb{E}[W^{\texttt{R}}]$, as in the standard GI/GI/1 model,  given the CDF $F_{W^{\texttt{R}}} $ of $W^{\texttt{R}}$ in  \eqref{eqn-FWR1}, one can calculate $\mathbb{E}[W^{\texttt{R}}] = \int_0^\infty (1-F_{W^{\texttt{R}}}(t) )dt$, but no concise and explicit formula is available. 

Although the above representation is relatively involved, it simplifies considerably in the special case where $X \geq 0$ almost surely. Specifically, applying \eqref{seriesRWX+} yields the following expression for the Laplace transform of $W^{\texttt{R}}$:
\begin{equation*}
\mathbb{E}\left[ \exp(- t W^{\texttt{R}} ) \right] = \frac{ q }{ 1- (1-q) \mathbb{E}\left[ \exp( -t X ) \right] }.
\end{equation*}
Letting $t \rightarrow \infty$, we obtain that
\begin{equation*}
F_{W^{\texttt{R}}} (0) = \frac{q}{1- (1-q) F_{X}(0)}.
\end{equation*}
Consequently,
\begin{equation*}
P_{\mathrm{delay}} = 1 - F_{W^{\texttt{R}}} (0) = \frac{ (1-q) ( 1- F_{X}(0) ) }{1- (1-q) F_{X}(0)}.
\end{equation*}
If $F_{X}(0)=0$, $P_{\mathrm{delay}} = 1-q$. 
Furthermore, when $X \geq 0$ almost surely, differentiating the Laplace transform yields explicit expressions for the mean and variance of the stationary waiting time $W^{\texttt{R}}$, under the assumptions $\mathbb{E}[X] < \infty$ and $\mathbb{E}[X^2] < \infty$, respectively. More specifically, 
\[
\mathbb{E}[W^{\texttt{R}}] = \frac{1-q}{q} \E[X]\,,
\]
\[
Var[W^{\texttt{R}}] =  \frac{1-q}{q} \E[X^2] +  \frac{(1-q)^2}{q^2} (\E[X])^2\,. 
\]
\end{remark}

 \vskip 0.8cm

\section{The GI/GI/$r$ queue with random resetting at arrival times}
\label{sec-GGr} 

\subsection{\it The Kiefer--Wolfowitz recursion for a {\rm{GI/GI/}}$r$ queue}\label{sec-GGr-standard}
\vskip .5cm
In a standard GI/GI/$r$ queue with $r>1$ servers and under the FCFS discipline, 
we operate with a collection of random vectors  $\{\unW_n\}$ where 
$\unW_n=(W_{n1},\ldots ,W_{nr})$ and $0\leq W_{n1}\leq\ldots\leq W_{nr}$.
In other words, $\unW_n$ takes values in the simplex $\bbS^+_{\leq}\subset\bbR^r$ where 
$\bbS^+_{\leq}=\Big\{\unx=(x_1,\ldots ,x_r):\;0\leq x_1\leq\ldots\leq x_r\Big\}$. 
Pictorially, $\unW_n$ represents the residual workload vector at the time
of arrival of the $n$th job, and its smallest entry, $W_{n1}$, gives the waiting time for the $n$th job.

The recursion that generates the sequence $\{\unW_n\}$ is due to Kiefer and Wolfowitz  \cite{KW1955}:
\beq
\label{eq:KWo1}\unW_{n+1}=\Big[\cR\big(\unW_n+V_n\,\une^{(1)}\big)-U_n\,{\underline{\,1\,}}\;\Big]^+
\eeq
with the following ingredients on the RHS: \vskip .5cm

\begin{itemize}
\item [(i)]
$V_n$ is the service time of the $n$th arrival, and $U_n$ is the time between the $n$th and $(n+1)$st
arrival.  It is assumed that the pairs $(U_n,V_n)$, $n=0,1,\ldots$, form an IID sequence.  The joint CDF
for $(U_n,V_n)$ is denoted by $G$: 
\beq\label{eqnG} G(u,v)=\mathbb{P}(U_n\leq u,V_n\leq v).\eeq 
We will assume that CDF $G$ is proper, i.e., RVs $U_n$ and $V_n$ take finite values only.
 \vskip .3cm

\item [(ii)]
$\une^{(1)}=(1,0,\ldots ,0)\in\bbZ_+^r$ and ${\underline{\,1\,}} =(1,\ldots ,1)\in\bbZ_+^r$ are $r$-dimensional 
$0,1$-vectors.\vskip .3cm

\item [(iii)]
$\cR\big(\unW_n+V_n\,\une^{(1)}\big)\in\bbS^+_{\leq}$ is the result of the re-arrangement operation $\cR$ 
applied to the
vector $\unW_n+V_n\,\une^{(1)}\in\bbR_+^r$: the vector $\cR\big(\unW_n+V_n\,\une^{(1)}\big)$ has
the same collection of entries as $\unW_n+V_n\,\une^{(1)}$ re-arranged in the non-decreasing order.
\vskip .3cm

\item [(iv)]
$\Big[\cR\big(\unW_n+V_n\,\une^{(1)}\big)-U_n\,{\underline{\,1\,}}\;\Big]^+\in\bbS^+_{\leq}$ is the vector 
obtained when the negative entries in $\cR\big(\unW_n+V_n\,\une^{(1)}\big)-U_k\,{\underline{\,1\,}}$ \ 
are replaced with zeros and non-negative entries are left intact. \vskip .3cm

\end{itemize}

Equation \eqref{eq:KWo1} generates a DTMC $\{\unW_n,\;n=0,1,\ldots\}$ on $\bbS^+_{\leq}$.
It can be re-written in terms of the $r$-dimensional CDFs 
$F_n(\unx )=\bbP(\unW_n\leq\unx )$, $n\geq 0$, as follows:
\begin{equation}
\label{eq:KWo2}
F_{n+1}(\unx )=\int_{\bbR^2}\int_{\bbS^+_{\leq}}{\bf 1}\Big(\unw\in\bbA(\unx ,u,v)\Big)\rd F_n(\unw )\rd G(u,v),\quad
\unx\in\RR^r,\end{equation} 
where 
the set  $\bbA(\unx ,u,v))\subset\bbS^+_{\leq}$ is given by
\begin{equation}
\bbA(\unx ,u,v)=\left\{\unw\in\bbS^+_{\leq}:\;\Big[\cR \big(\unw +v\une^{(1)}\big)-u\,{\underline{\,1\,}}\;\Big]^+\leq\unx\right\},
\end{equation} 
and $G(u,v)$ is given in \eqref{eqnG}. Here and below, the inequality  between vectors means the inequality 
between their respective entries.

As before, it is instructive to write equation \eqref{eq:KWo2} in an operator form:
\beq\label{OpKfr}\beal
F_{n+1}={\tt K}F_n\quad\hbox{where operator ${\tt K}$ acts on a CDF $H$ by}\\
\qquad ({\tt K}H)(\unx )=\diy\int_{\bbR^2_+}\int_{\bbS^+_{\leq}}{\bf 1}\Big(\unw\in\bbA (\unx ,u,v)\Big)\rd H(\unw )
\rd G(u,v),\quad\unx\in\RR^r.\ena\eeq

It is known that if the  traffic intensity $\rho :=\dfrac{\bbE [V]}{r\bbE [U]}<1$,  then the 
stationary Kiefer--Wolfowitz equation 
\beq\label{eq:KWoS1}\beal
\unW \overset{d}{=} \Big[\cR\big(\unW+V\,\une^{(1)}\big)-U\,\underline{\, 1 \,}\;\Big]_+
\quad\hbox{or, equivalently, $F={\tt K}F$, \hbox{i.e.,}}\\
\hskip 2cm F(\unx )=\diy\int_{\bbR^2}\int_{\bbS^r}{\bf 1}(\unw\in A(\unx ,u,v))
\rd F(\unw )\rd G(u,v)\ena\eeq 
has a unique solution giving a proper CDF $F$ on $\bbR^r$. In fact, for $\rho <1$, the DTMC $\{\unW_n\}$ is Harris 
ergodic. On the other hand, when $\rho\geq 1$, there is no proper CDF $F$ on $\RR^r$ satisfying 
equation \eqref{eq:KWoS1}, Cf.  \cite[Chapter XII.2]{asmussen2003}.\vskip .5cm

\subsection{\it The modified Kiefer--Wolfowitz recursion for a {\rm{GI/GI/}}$r$ queue with resetting.}\label{sec-GGr-reset}
\vskip .5cm
The model with random resetting at arrival times again involves the parameter $q\in (0,1)$. Set 
${\underline{\,0\,}}=(0,\ldots ,0)$. Equations \eqref{eq:KWo1} and
 \eqref{eq:KWo2} are replaced with
\beq\label{eq:KWoR1}\unW^\tR_{n+1}=\bcs \qquad\qquad{\underline{\,0\,}},&\hbox{with probability $q$,}\\
\Big[\cR\big(\unW^\tR_n+V_n\,\une^{(1)}\big)-U_n\,{\underline{\,1\,}}\;\Big]^+,&\hbox{with probability $1-q$,}\ecs\eeq
and
\beq\label{eq:KWoR2}\beal F^\tR_{n+1}(\unx )=q{\bf 1}_{\bbS^+_{\leq}}(\unx ) +(1-q )\diy\int_{\bbR^2}\int_{\bbS^+_{\leq}}{\bf 1}\Big(\unw\in\bbA(\unx ,u,v)\Big)\rd F^\tR_n(\unw )
\rd G(u,v),\ena\eeq
respectively, with $F^\tR_n(\unx )=\bbP (\unW_n\leq\unx)$. As above, equation \eqref{eq:KWoR1} determines 
a DTMC $\{\unW^\tR_n,\;n=0,1,\ldots\}$ on 
$\bbS^+_{\leq}$. Again, the vector $\unW^\tR_n$ represents the residual workloads at the servers at the  
$n$th arrival time. 

Accordingly, the stationary equations for $\{\unW^\tR_n\}$ take the following equivalent forms:
\beq\label{eq:KWoRS}\beal\;\;{\rm{(i)}}\quad \unW^\tR \overset{d}{=} \bcs \qquad\qquad{\underline{\,0\,}},&\hbox{with 
probability $q$,}\\
\Big[\cR\big(\unW^\tR+V\,\une^{(1)}\big)-U\,{\underline{\,1\,}}\;\Big]^+,&\hbox{with probability $1-q$,}\ecs\quad\hbox{or}\\
\beal{\rm{(ii)}}\quad F^\tR(\unx )=q{\bf 1}_{\bbS^+_{\leq}}(\unx )\\
\hskip 2cm +(1-q )\diy\int_{\bbR^2}\int_{\bbS^+_{\leq}}{\bf 1}\Big(\unw\in\bbA(\unx ,u,v)\Big)\rd F^{\tR} 
(\unw )\rd G(u,v),\ena\;\;\hbox{or}\\
{\rm{(iii)}}\quad F^\tR=q{\bf 1}_{\bbS^+_{\leq}}+(1-q){\tt K}F^\tR\quad\Llra\quad
\big({\tt I}-(1-q){\tt K}\big)F^\tR=q{\bf 1}_{\bbS^+_{\leq}}\\
\hskip 1.2cm \hbox{solved by}\\
{\rm{(iv)}}\quad F^\tR=q\big({\tt I}-(1-q){\tt K}\big)^{-1}{\bf 1}_{\bbS^+_{\leq}}\\
\hskip 4 cm =q{\bf 1}_{\bbS^+_{\leq}}
+q(1-q){\tt K}{\bf 1}_{\bbS^+_{\leq}}+q(1-q)^2{\tt K}^2{\bf 1}_{\bbS^+_{\leq}}+\ldots .\ena\eeq
\vskip .5cm

\begin{prop}\label{Prop6.1}
For any $q\in (0,1)$ and a sequence of \ {\rm{IID}} \ {\rm{RV}} pairs \ $\{(U_n,V_n)\}$,  
the {\rm{DTMC}} \ $\{\unW^\tR_{n}\}$ \ is Harris ergodic. 
\end{prop}

\begin{proof} We will only give here a sketch of the (rather tedious) proof as it does not contain
serious novel elements.  It is reduced to a repetition of arguments from 
\cite[Chapters XII.1, XII.2]{asmussen2003}. The crux of the matter is Theorem 1.2 on page 432 in
\cite[Chapter XII.1]{asmussen2003} and Theorem 2.2 on page 345 in \cite[Chapter XII.2]{asmussen2003}
rewritten in a modified form for the DTMC with resetting $\{W^\tR_{n}\}$. In turn, the proof of the modified 
theorems is based on analogs of Lemma 1.3 and Lemmas 2.3 and 2.4 in \cite[Chapters XII.1, XII.2]{asmussen2003}.
Such analogs connect the DTMC $\{W^\tR_{n}\}$ with the majorizing Markov chain $\{{\wt W}^\tR_{n}\}$ where 
arriving jobs are directed to servers in the cyclic order with probability $1-q$ and trigger resetting of
the whole vector of residual workloads to $\underline{\,0\,}$ with probability $q$. The analysis of the  majorizing DTMC  $\{{\wt W}^\tR_{n}\}$ is essentially reduced to the GI/GI/$1$ model with resetting which leads to the
assertion of Proposition \ref{Prop6.1}.
\end{proof}

The above construction then leads to the following result. 

\begin{theorem} \label{thm-GGr}
For any $q\in (0,1)$, the series in \eqref{eq:KWoRS}{\rm{(iv)}} determines 
a proper {\rm{CDF}} satisfying \eqref{eq:KWoRS}{\rm{(ii)}}.
Furthermore, equation \eqref{eq:KWoRS}{\rm{(ii)}} has a unique bounded solution, and this solution 
is given by the series in \eqref{eq:KWoRS}{\rm{(iv)}}.
\end{theorem}



\begin{proof}
As in the case of the model GI/GI/$1$ with resetting, the fact that the series in \eqref{eq:KWoRS}{\rm{(iv)}}
gives a solution to \eqref{eq:KWoRS}{\rm{(ii)}} follows from the construction with the help of  
Proposition \ref{Prop6.1}. Uniqueness is also established by the same argument as for the GI/GI/$1$
model. 
\end{proof}


\medskip
\section{The GI/GI/$\infty$ queue with random resetting at arrival times}
 \label{sec-ISQ}
 
 \medskip

\subsection{\it The recursion for a {\rm{GI/GI/}}$\infty$ queue.}  \label{rec-ISQ-arrival}

In the section, we study GI/GI/$\infty$ queue, where infinitely many servers are available.
A standard GI/GI/$\infty$ can be described via a sequence of random-dimension random vectors 
$\unW_n=(W_{n1},\ldots ,W_{nS(n)})$, $n=0,1,\dots$. Specifically, $S(n)$ is a random variable with nonnegative integer values, which determines the dimension of $\unW_n$ and corresponds to the number of active servers when the $n$th job arrives. If $S(n) \geq 1$, the entries of $\unW_n$ satisfy that $W_{n1}\geq\ldots\geq W_{nS(n)}>0$, and they correspond to the residual workloads of the active servers. In particular, the largest entry, $W_{n1}$, gives the time needed for clearing the system of jobs entered before the $n$th arrival time. If $S(n)=0$, then $\unW_n$ reduces to $\mathfrak{E}$, meaning the $n$th job finds an empty queue at its arrival. With the above definition, it is straightforward to see that $\unW_n$ takes values in the union 
${\un\bbO}_+(\geq ):=\onwl{\cup}\limits_{s= 0}^{\infty} \bbO^s_+(\geq )$ of simplexes 
$\bbO^s_+(\geq )=\big\{\unx=(x_1,\ldots ,x_s):\;x_1\geq\ldots\geq x_s>0\big\}$ of varying 
dimension $s\geq 1$, and a single-state set $\bbO^0_+(\geq )=\{\mathfrak{E}\}$ for $s=0$.

The recursion for the residual workload vector $\{\unW_n\}$  in the GI/GI/$\infty$ model is given by
\beq\label{eq:recISQ}\unW_{n+1}=\cR\Big(\cS\left\{\big[\cP(V_n,\unW_n)-U_n{\underline{\,1\,}}_{S(n)+1}\big]^+\right\}\Big),\quad n=0,1,\ldots .\eeq
Here the RHS contains the following components: \vskip .5cm

\begin{itemize}
\item [(i)]
As in Section \ref{sec-GGr}, $V_n$ is the service time of the $n$th arrival, and $U_n$ is the time between the $n$th 
and $(n+1)$st arrival. It is assumed that the pairs $(U_n,V_n)$, $n=0,1,\ldots$, form an IID sequence.  The 
joint CDF for $(U_n,V_n)$ is again denoted by $G$ and assumed to be proper.  
 \vskip .3cm

\item [(ii)] $S(n)$ is the dimension of $\unW_n$ and the vector ${\underline{\,1\,}}_{S(n)+1}=(1,\dots ,1)\in\bbZ_+^{S(n)+1}$ has all entries $1$. \vskip .3cm

\item [(iii)] $\cP(V_n,\unW_n)\in\bbR^{S(n)+1}$ is the 
result of concatenation of the value $V_n$ and the vector $\unW_n$. \vskip .3cm

\item [(iv)] $\big[\cP(V_n,\unW_n)-U_n{\underline{\,1\,}}_{S(n)+1}\big]^+$ is the vector obtained from 
$\cP(V_n,\unW_n)-U_n{\underline{\,1\,}}_{S(n)+1}$
when negative entries are replaced with zeros and non-negative entries are left intact. \vskip .3cm

\item [(v)] $\cS\left\{\big[\cP(V_n,\unW_n)-U_n{\underline{\,1\,}}_{S(n)+1}\big]^+\right\}$ is the result of shortening 
vector $\big[\cP(V_n,\unW_n)-U_n{\underline{\,1\,}}_{S(n)+1}\big]^+$ by removing the zero entries. If all entries are removed, it is denoted by $\mathfrak{E}$.
Furthermore, the dimension of $\cS\Big\{\big[\cP(V_n,\unW_n)-U_n{\underline{\,1\,}}_{S(n)+1}\big]^+\Big\}$ equals
$S(n+1)$. \vskip .3cm 

\item [(iv)] $\cR\Big(\cS\left\{\big[\cP(V_n,\unW_n)-U_n{\underline{\,1\,}}_{S(n)+1}\big]_+\right\}
\Big)$ is the result of the re-arrangement  applied to the
vector $\cS\left\{\big[\cP(V_n,\unW_n)-U_n{\underline{\,1\,}}_{S(n)+1}\big]^+\right\}$.
Indeed, it
has the same collection of entries as $\cS\left\{\big[\cP(V_n,\unW_n)-U_n{\underline{\,1\,}}_{S(n)+1}\big]^+\right\}$ re-arranged in the non-increasing order.
\vskip .3cm
\end{itemize}

Equation \eqref{eq:recISQ} generates a DTMC $\{\unW_n,\;n=0,1,\ldots\}$ on ${\un\bbO}_+(\geq )$.
The probability distribution of $\unW_n$ on ${\un\bbO}_+(\geq )$ is described by a sequence 
$\unF_n =(F_{n,0},F_{n,1},F_{n,2},\ldots )$ 
where
\beq F_{n,0}=\bbP\big(\unW_n= \mathfrak{E} \big), \eeq
and for each $k \ge 1$, $F_{n,k}:  \bbO^k_+(\geq ) \to [0,1]$ is a nondecreasing, right-continuous function
such that 
\beq F_{n,k}(\unx^{(k)})=\bbP\big(S(n)=k,\unW_n\leq\unx^{(k)}\big),
\;\;\unx^{(k)}\in  \bbO^k_+(\geq ).\eeq
Let $ F_{n,k}^* := \sup_{\unx^{(k)} \in \bbO^k_+(\geq ) }  F_{n,k}(\unx^{(k)})$ for each $k\ge 1$. Then it must hold that 
\[
F_{n,0} + F_{n,1}^* + \cdots + F_{n,k}^* + \cdots =1. 
\]


Equation \eqref{eq:recISQ} can be re-written in terms of the sequences
$\unF_n$, as follows:
\begin{equation}
\label{eq:KWo2-ISQ}\beal
F_{n+1,0}=\sum\limits_{l=0}^{\infty}\diy\int_{\bbR^2_+}\int_{\bbO^l_+(\geq )}
{\bf 1}\Big(\unw\in\bbA_{0,l}(u,v)\Big)\rd F_{n,l}(\unw )\rd G(u,v),\\
F_{n+1,k}(\unx^{(k)})=\sum\limits_{l=k-1}^{\infty}\diy\int_{\bbR^2_+}\int_{\bbO^l_+(\geq )}
{\bf 1}\Big(\unw\in\bbA_{k,l}(\unx^{(k)} ,u,v)\Big)\rd F_{n,l}(\unw )\rd G(u,v),\\
\hskip 7.5 cm\unx^{(k)}\in  \bbO^k_+(\geq ) ,\;\;k\geq 1.\ena\end{equation} 
Here, the sets  $\bbA_{0,l}(u,v),\bbA_{k,l}(\unx^{(k)},u,v))\subset\bbO^l_+(\geq )$ are given by
\begin{equation}\beal
\bbA_{0,l}(u,v)=\bigg\{\unw^{(l)}\in\bbO^l_+(\geq ): \cR\Big(\cS\left\{\big[\cP(v,\unw^{(l)})
-u{\underline{\,1\,}}_{l+1}\big]^+\right\}\Big)=\mathfrak{E} \bigg\},\\
\bbA_{k,l}(\unx^{(k)},u,v)=\bigg\{\unw^{(l)}\in\bbO^l_+(\geq ): \cR\Big(\cS\left\{\big[\cP(v,\unw^{(l)})
-u{\underline{\,1\,}}_{l+1}\big]^+\right\}\Big)\in\bbO^k_+(\geq ),\\
\hskip 5cm\cR\Big(\cS\left\{\big[\cP(v,\unw^{(l)})
-u{\underline{\,1\,}}_{l+1}\big]^+\right\}\Big) \leq\unx^{(k)}\bigg\},\quad k\geq 1.
\ena\end{equation} 

As before, it is instructive to write equation \eqref{eq:KWo2-ISQ} in an operator form:
\beq
F_{n+1,k}=\sum\limits_l{\tt K}_{k,l}F_{n,l},
\eeq
where operator ${\tt K}_{k,l}$ acts on a nondecreasing right continuous function $H_l: \bbO^l_+(\geq ) \to [0,1]$  by 
\beq
\beal
{\tt K}_{0,l}H_l =\diy\int_{\bbR^2_+}\int_{\bbO^l_+(\geq )}{\bf 1}\Big(\unw^{(l)}\in\bbA_{0,l} 
(u,v)\Big)\rd H(\unw^{(l)})\rd G(u,v),\\
({\tt K}_{k,l}H_l)(\unx^{(k)})=\diy\int_{\bbR^2_+}\int_{\bbO^l_+(\geq )}{\bf 1}\Big(\unw^{(l)}\in\bbA_{k,l} 
(\unx^{(k)},u,v)\Big)\rd H(\unw^{(l)})\rd G(u,v),\\
\hskip 7cm\unx^{(k)}\in  \bbO^k_+(\geq ) ,\;k\geq 1.\ena
\eeq


For future use, it is convenient to introduce the operator $\un{\,{\tt K}\,}=({\tt K}_{k,l})$ with blocks ${\tt K}_{k,l}$ acting on the sequences $\unH=(H_0,H_1,H_2,\ldots )$:
\beq\beal \un{\,{\tt K}\,}\,\unH =\Big(\left(\un{\,{\tt K}\,}\,\unH\right)_0,\left(\un{\,{\tt K}\,}\,\unH\right)_1,
\left(\un{\,{\tt K}\,}\,\unH\right)_2,\ldots \Big)\\
\qquad\hbox{where}\quad
\left(\un{\,{\tt K}\,}\,\unH\right)_k=\sum\limits_{l=[k-1]^{+}}^{\infty}{\tt K}_{k,l}H_l,\quad k=0,1,\ldots .
\ena\eeq

\subsection{\it The recursion for a {\rm{GI/GI/}}$\infty$ queue with resetting.} 
\label{res-ISQ-arrival}

The recursion for a GI/GI/$\infty$ model with resetting takes the form
\beq\label{rec:ISQres}\unW^\tR_{n+1}=\bcs \qquad \mathfrak{E},&\hbox{with probability $q$,}\\
\cR\Big(\cS\left\{\big[\cP(V_n,\unW_n)-U_n{\underline{\,1\,}}_{S(n)+1}\big]^+\right\}\Big),
&\hbox{with probability $1-q$.}\ecs\eeq
It generates a DTMC $\{\unW^\tR_n,\;n=0,1,\ldots\}$ on ${\un\bbO}_+(\geq )$.

As above, we rewrite equation \eqref{rec:ISQres} in terms of the sequence of functions  $\unF^\tR_n =(F^\tR_{n,0},F^\tR_{n,1},\ldots )$, 
where 
$$
F^\tR_{n,0}  = \bbP\Big(\unW^\tR_n = \mathfrak{E}\Big),
$$
and for $k\geq 1$, $F^\tR_{n,k}: \bbO^k_+(\geq ) \to [0,1]$ is a nondecreasing, right-continuous function
such that 
$$F^\tR_{n,k} ( \unx^{(k)} ) =\bbP\Big(\unW^\tR_n\in\bbO^k_+(\geq ),\;\unW_n\leq\unx^{(k)}\Big).$$ 
Denote by $ \underline{\,e\,}^{(1)} = (1,0,0, \dots) $ the infinite-dimensional vector whose first entry is $1$ and all remaining entries are $0$.
Then, we have
\begin{align*}
\unF^\tR_{n+1}&=q\,\underline{\,e\,}^{(1)} +(1-q)\,{\un{\,\tK\;}}\,{\unF}^\tR_n,
\quad\hbox{or, entry-wise,}\\
F^{\texttt{R}}_{n+1, 0} &= q  + (1-q) \sum_{l=0}^{\infty}  \tK_{0,l}F^\tR_{n,l}, \\
F^{\texttt{R}}_{n+1, k} ( \underline{x}^{(k)} ) &= (1-q) \sum_{l=k-1}^{\infty}  (\tK_{k,l}F^\tR_{n,l})(\unx^{(k)}),  \,\unx^{(k)}\in \bbO^k_+(\geq ), \quad k=1, 2,\dots.
\end{align*}
Therefore, the stationary version becomes 
\beq\label{resIQSeq}\beal 
\unF^\tR=q\,{\underline{\,e\,}}^{(1)} +(1-q)\,{\un{\,\tK\;}}\,{\unF}^\tR
\Llra ({\un{\,{\tt I}\,}}-(1-q){\un{\,{\tt K}\,}}){\unF}^\tR =q\,{\underline{\,e\,}}^{(1)}.\ena\eeq
Equation \eqref{resIQSeq} can be solved by 
\beq\label{eq7.10}\unF^\tR=q\Big({\underline{\,e\,}}^{(1)}
+(1-q){\un{\,{\tt K}\,}} \, {\underline{\,e\,}}^{(1)}
+(1-q)^2{\un{\,{\tt K}\,}}^2
{\underline{\,e\,}}^{(1)}+\cdots \Big).
\eeq 

Therefore we conclude the following result. 

\begin{theorem}
The DTMC $\{\unW^\tR_n\}$ has a stationary probability distribution characterized by equation \eqref{eq7.10}.
\end{theorem}

\medskip

\section{Concluding Remarks}
\medskip

In this paper,  we have considered the standard queueing models with random resettings.
Several extensions are possible future works.
First, an immediate extension would be to consider more general Markov chains with random resettings.
 It would be interesting to identify conditions under which an explicit stationary distribution could be derived. 
Some efforts in this direction are made in our forthcoming paper  \cite{pang2025Wiener}.
Second, for non--Markovian queues, it would be interesting to consider different forms of resettings other than those at arrival times.  
Third, diffusions have been established to approximate the performances of queues in heavy traffic.  
It would be also worth considering such diffusion models with random resetting, particularly, their ergodic properties and characterization of stationary distributions.

\medskip

{\bf Acknowledgements.} G. Pang is partly supported by NSF grants DMS 2216765 and CMMI 2452849. 
 I. Stuhl and Y. Suhov thank Math Dept, Penn State University for
support. Y. Suhov thanks DPMMS, University of Cambridge, and St John's College, Cambridge,
for support. Y. Suhov thanks IHES, Bures-sur-Yvette, for hospitality during a visit in 2024.

\bigskip


\section{Appendix: Additional Proofs for the M/M/$r$ Queue with Resetting} \label{sec-appendix}

The appendix is devoted to the proofs of Lemmas \ref{lem-pi-positive}--\ref{lm:interationAC}. Since the proof of Lemma~\ref{lem-pi-positive} relies on the results established in Lemmas~\ref{lm:recursionforAn}--\ref{lm:interationAC}, we first prove Lemmas~\ref{lm:recursionforAn}--\ref{lm:interationAC} and then use these results to establish Lemma~\ref{lem-pi-positive}.

We begin with the proof of Lemma~\ref{lm:recursionforAn}.

\begin{proof}[Proof of Lemma~\ref{lm:recursionforAn}]
By integration by parts, we obtain
\begin{align*}
A_0 &= \int_0^1 \gamma (1-s)^{\gamma-1} e^{-\theta s}ds = \int_0^1 e^{-\theta s} d (- (1-s)^\gamma) \\ 
&= 1+ \int_0^1(1-s)^\gamma d \left(e^{-\theta s}\right) =  1-\theta \int_0^1(1-s)^\gamma  e^{-\theta s}ds  \\
&= 1- \theta  \int_0^1(1-s)^{\gamma -1}  e^{-\theta s}ds + \theta  \int_0^1(1-s)^{\gamma -1} s  e^{-\theta s}ds \\
& = 1 - \frac{\theta}{\gamma} A_0 + \frac{\theta}{\gamma} A_1,
\end{align*}
which immediately implies \eqref{eqinf:A0A1}. We now turn to \eqref{eqn-An-recursion}. Again, applying integration by parts, it follows
\begin{align*}
A_{n} &= \int_0^1 \gamma (1-s)^{\gamma-1} \theta s^n e^{-\theta s}ds   
= \int_0^1 \theta s^n e^{-\theta s}d (- (1-s)^\gamma) \\
&=   \int_0^1 (1-s)^\gamma  \theta  ( n s^{n-1} e^{-\theta s} - s^n \theta e^{-\theta s})ds \\
&= n \int_0^1 (1-s)^\gamma  \theta s^{n-1} e^{-\theta s} \, ds - \theta \int_0^1 (1-s)^\gamma  \theta s^{n} e^{-\theta s} \, ds \\
& = n \left( \int_0^1 (1-s)^{\gamma -1}  \theta s^{n-1} e^{-\theta s} \, ds - \int_0^1 (1-s)^{\gamma -1}  \theta s^{n} e^{-\theta s} \, ds  \right) \\
& \quad - \theta\left(\int_0^1 (1-s)^{\gamma -1}  \theta s^{n} e^{-\theta s} \, ds -\int_0^1 (1-s)^{\gamma -1}  \theta s^{n+1} e^{-\theta s} \, ds  \right) \\
&= n\left( \frac{1}{\gamma} A_{n-1} - \frac{1}{\gamma} A_{n} \right) - \theta \left( \frac{1}{\gamma} A_{n}  - \frac{1}{\gamma} A_{n+1} \right).
\end{align*}
This is equivalent to  the recursive equation in \eqref{eqn-An-recursion}.
\end{proof}

Subsequently, we prove Lemma~\ref{lm:recursionforCk}.

\begin{proof}[Proof of Lemma~\ref{lm:recursionforCk}]
Since \eqref{eqinf:C0C1} can be verified by direct calculation, it suffices to prove the recursive formula \eqref{eqn-Ck-recursion} and \eqref{eqr:thedenominator>0}.

We begin with \eqref{eqn-Ck-recursion}.
By recalling the definition of $C_{k}$ in \eqref{eqn-Ck-MMr} and applying the identity
\begin{equation*}
\binom{k+1}{l} \mathds{1}_{\{l=0,1,\dots, k+1\}}= \binom{k}{l} \mathds{1}_{\{l=0,1,\dots, k\}} + \binom{k}{l-1} \mathds{1}_{\{l=1,2,\dots, k+1\}},
\end{equation*}
it follows that
\begin{eqnarray*}
\theta C_{k+1} &=& \theta  \sum_{l=0}^{k+1} \binom{k+1}{l} \gamma (\gamma +1) \cdots (\gamma + l-1) \theta^{-l} \\
&=& \theta  \sum_{l=0}^{k} \binom{k}{l} \gamma (\gamma +1) \cdots (\gamma + l-1) \theta^{-l} \\
&& + \theta  \sum_{l=1}^{k+1} \binom{k}{l-1} \gamma (\gamma +1) \cdots (\gamma + l-1) \theta^{-l} \\
&=& \theta C_{k} + \theta  \sum_{l=1}^{k+1} \binom{k}{l-1} \gamma (\gamma +1) \cdots (\gamma + l-1) \theta^{-l}.
\end{eqnarray*}
Furthermore, we have
\begin{eqnarray*}
&& \theta  \sum_{l=1}^{k+1} \binom{k}{l-1} \gamma (\gamma +1) \cdots (\gamma + l-1) \theta^{-l}  \\
&=&  \theta  \sum_{l=0}^{k} \binom{k}{l} \gamma (\gamma +1) \cdots (\gamma + l) \theta^{-l-1} \\
&=&    \sum_{l=0}^{k} \binom{k}{l} \gamma (\gamma +1) \cdots (\gamma + l -1) (\gamma + l) \theta^{-l} \\
&=& \sum_{l=0}^{k} \binom{k}{l} \gamma (\gamma +1) \cdots (\gamma + l -1)  \times \gamma  \theta^{-l} \\
&& + \sum_{l=0}^{k} \binom{k}{l} \gamma (\gamma +1) \cdots (\gamma + l -1)  \times l \theta^{-l} \\
&=&  \gamma  \sum_{l=0}^{k} \binom{k}{l} \gamma (\gamma +1) \cdots (\gamma + l -1)  \theta^{-l} \\
&& +  \sum_{l=0}^{k} l\binom{k}{l} \gamma (\gamma +1) \cdots (\gamma + l -1)  \theta^{-l} \\
&=&  \gamma C_{k} +  \sum_{l=0}^{k} l\binom{k}{l} \gamma (\gamma +1) \cdots (\gamma + l -1)  \theta^{-l},
\end{eqnarray*}
where in the first equality, we substitute $l-1$ with $l$. Combining the last two displays yields 
\begin{align*}
\theta C_{k+1} = (\theta +\gamma) C_{k} + \sum_{l=0}^{k} l\binom{k}{l} \gamma (\gamma +1) \cdots (\gamma + l -1)  \theta^{-l}.
\end{align*}
Comparing the last display with the recursive formula \eqref{eqn-Ck-recursion} for $C_{k}$, it remains to prove
\begin{equation*}
k C_{k} = k C_{k-1} + \sum_{l=0}^{k} l\binom{k}{l} \gamma (\gamma +1) \cdots (\gamma + l -1)  \theta^{-l}.
\end{equation*}
Observing that the term inside the summation is zero when $l=0$ and by applying the identity $ l \binom{k}{l} = k \binom{k-1}{l-1} $, we have
\begin{eqnarray*}
&& \sum_{l=0}^{k} l\binom{k}{l} \gamma (\gamma +1) \cdots (\gamma + l -1)  \theta^{-l} 
= \sum_{l=1}^{k} l\binom{k}{l} \gamma (\gamma +1) \cdots (\gamma + l -1)  \theta^{-l} \\
&=& \sum_{l=1}^{k} k \binom{k-1}{l-1} \gamma (\gamma +1) \cdots (\gamma + l -1)  \theta^{-l}
= k \sum_{l=1}^{k} \binom{k-1}{l-1} \gamma (\gamma +1) \cdots (\gamma + l -1)  \theta^{-l}.
\end{eqnarray*}
Therefore,
\begin{eqnarray*}
&&  k C_{k-1} + \sum_{l=0}^{k} l\binom{k}{l} \gamma (\gamma +1) \cdots (\gamma + l -1)  \theta^{-l} \\
&=& k \sum_{l=0}^{k-1} \binom{k-1}{l} \gamma (\gamma +1) \cdots (\gamma + l -1)  \theta^{-l} \\
&& + k \sum_{l=1}^{k} \binom{k-1}{l-1} \gamma (\gamma +1) \cdots (\gamma + l -1)  \theta^{-l} \\
&=& k \sum_{l=0}^{k} \binom{k}{l} \gamma (\gamma +1) \cdots (\gamma + l -1)  \theta^{-l} \\
&=& k C_{k},
\end{eqnarray*}
where the second last equality follows by the identity 
\begin{equation*}
\binom{k-1}{l-1} \mathds{1}_{\{l=1,2,\dots, k\}} + \binom{k-1}{l}\mathds{1}_{\{l=0,1,\dots, k-1\}} = \binom{k}{l} \mathds{1}_{\{l=0,1,\dots, k\}}.
\end{equation*}
We complete th proof of the recursive formula \eqref{eqn-Ck-recursion}.

We then turn to \eqref{eqr:thedenominator>0}. When $\theta >0$, it is straightforward to verify that $\theta > r \alpha >0$ and
\begin{align*}
C_{r} &= \sum_{l=0}^{r} \binom{r}{l} \gamma (\gamma +1) 	\cdots (\gamma + l -1) \theta^{-l}  \\
&= \sum_{l=0}^{r-1} \binom{r-1}{l} \gamma (\gamma +1) 	\cdots (\gamma + l -1) \theta^{-l}  + \sum_{l=1}^{r} \binom{r-1}{l-1} \gamma (\gamma +1) 	\cdots (\gamma + l -1) \theta^{-l} \\
&> \sum_{l=0}^{r-1} \binom{r-1}{l} \gamma (\gamma +1) 	\cdots (\gamma + l -1) \theta^{-l} \\
&= C_{r-1} 
 >0,
\end{align*} 
where in the first equality, we apply the following identity:
\begin{equation*}
\binom{r}{l} \mathds{1}_{\{ l=0,1, \dots, r \}} = \binom{r-1}{l} \mathds{1}_{\{ l=0,1, \dots, r-1 \}} + \binom{r-1}{l-1} \mathds{1}_{\{ l=1, 2, \dots, r \}}.
\end{equation*}
Therefore, \eqref{eqr:thedenominator>0} follows directly.
\end{proof}

We then turn to the proof of Lemma~\ref{lm:interationAC}.

\begin{proof}[Proof of Lemma~\ref{lm:interationAC}]
By applying the recursive formulas for $A_k$ in \eqref{eqn-An-recursion} and $C_{k}$ in \eqref{eqn-Ck-recursion}, we obtain
\begin{eqnarray*}
&&  \theta^{r} (  A_{r-1} C_{r} -  C_{r-1} A_{r}) \\
&=& \theta^{r-1} (  A_{r-1} \theta C_{r} -  C_{r-1}\theta A_{r}) \\
&=& \theta^{r-1}  A_{r-1} \left[  (r-1+\theta+\gamma) C_{r-1} - (r-1) C_{r-2} \right]  \\
&& -\theta^{r-1} C_{r-1} \left[  (r-1+\theta+\gamma) A_{r-1} - (r-1) A_{r-2} \right]  \\
&=& (r-1) \theta^{r-1}  (  A_{r-2} C_{r-1} -  C_{r-2} A_{r-1}).
\end{eqnarray*}
Furthermore, It follows from \eqref{eqinf:A0A1} and \eqref{eqinf:C0C1} that $\theta(A_0  C_1 - C_0 A_1) = \gamma$. By combining these two results, we directly obtain that $ \theta^{r} (  A_{r-1} C_{r} -  C_{r-1} A_{r}) = \gamma (r-1)! $.
\end{proof}

Finally, we prove Lemma~\ref{lem-pi-positive}.

\begin{proof}[Proof of Lemma \ref{lem-pi-positive}]
We first show that the sequence $\pi_k$ defined in \eqref{eqn-pi-MMr-1} and \eqref{eqn-pi-MMr-2} is convergent. Note that $\pi_{k}  = \alpha^{k-r+1} \pi_{r-1} $ for $k=r,r+1, \dots$. It is immediate to check that $0 < \alpha <1$ when $\theta ,\gamma >0$. Therefore, to prove the convergence, it suffices to prove that $\pi_0, \pi_1, \dots, \pi_{r-1}$ are finite. However, it follows from \eqref{eqr:thedenominator>0} that $L_{r-1,r}$ is finite. Together with the fact that $A_{k}$ and $C_{k}$ are finite for $k=0,1,\dots$, we obtain the finiteness of $\pi_0, \pi_1, \dots, \pi_{r-1}$.



Next, we prove that the $\pi_k$ defined in \eqref{eqn-pi-MMr-1} and \eqref{eqn-pi-MMr-2} are always nonnegative. 
It suffices to show that
\begin{equation*}
A_{k } -L_{r-1 ,r} C_{k} \geq 0 \qquad \text{for $k=0,1,\dots, r-1$}.
\end{equation*}
To simplify the notation, we define $D_{k} := A_{k } -L_{r-1 ,r} C_{k}$ for $k=0,1,\dots$. 
Then, by the expression of $L_{r-1 ,r}$ in \eqref{eqn-L-MMr}, we obtain
\begin{equation*}
D_{k}= A_{k} - \frac{  \theta A_{r} - r \alpha A_{r-1} }{ \theta C_{r} - r \alpha C_{r-1} } C_{k}
= \frac{ \theta A_{k} C_{r} - \theta C_{k} A_{r} - r\alpha A_{k}C_{r-1} +  r\alpha C_{k}A_{r-1} }{\theta C_{r} - r \alpha C_{r-1}}.
\end{equation*}
In particular, when $k =r$,
\begin{equation*}
D_{r} = \frac{ r \alpha ( A_{r-1} C_{r} - C_{r-1} A_{r} ) }{ \theta C_{r} - r \alpha C_{r-1} },
\end{equation*}
and when $k = r-1$,
\begin{equation*}
D_{r-1} = \frac{ \theta ( A_{r-1} C_{r} - C_{r-1} A_{r} ) }{ \theta C_{r} - r \alpha C_{r-1} }.
\end{equation*}
Then, the positivity of $D_{r}$ and $D_{r-1}$ follows from \eqref{eqr:keyidentityAC} and \eqref{eqr:thedenominator>0}. Indeed, we have $D_{r-1} > D_{r} >0$, since $\theta > r \alpha$. Using the recursive formulas for $A_k$ and $C_k$, we obtain that $D_{k}$ satisfies the same recursive formula, that is,
\begin{equation*}
\theta D_{k+1 } - (k + \theta + \gamma)   D_k +  k  D_{k-1} =0, \qquad \text{for $k =1,2, \dots$.}
\end{equation*}
Then, by mathematical induction, we can show that $D_{k-1} > D_{k} >0 $ for $k=1,2, \dots, r$. Indeed, if $D_{n-1} > D_{n} >0 $, then
\begin{align*}
(n-1) ( D_{n-2} - D_{n-1} ) = \theta (D_{n-1} - D_{n} ) + \gamma D_{n-1} >0,
\end{align*}
which gives $D_{n-2} > D_{n-1} >0$. 
Thus, we conclude the nonnegativity of $\pi_k$.
\end{proof}

\end{document}